\documentclass[final,onefignum,onetabnum]{siamonline220329}

\usepackage{lipsum}
\usepackage{amsfonts}
\usepackage{graphicx}
\usepackage{epstopdf}
\usepackage{algorithmic}
\ifpdf
  \DeclareGraphicsExtensions{.eps,.pdf,.png,.jpg}
\else
  \DeclareGraphicsExtensions{.eps}
\fi

\usepackage{enumitem}
\setlist[enumerate]{leftmargin=.5in}
\setlist[itemize]{leftmargin=.5in}

\newsiamremark{remark}{Remark}
\newsiamremark{hypothesis}{Hypothesis}
\crefname{hypothesis}{Hypothesis}{Hypotheses}
\newsiamthm{claim}{Claim}

\headers{Semidefinite Relaxations in Robust Rotation Search}{L. Peng, M. Fazlyab, and R. Vidal}

\title{Towards Understanding The Semidefinite Relaxations of Truncated Least-Squares in Robust Rotation Search\thanks{This work is an extension of our conference paper \cite{Peng-ECCV2022}.
		\funding{This work was supported by grants NSF 1704458, NSF 1934979 and ONR MURI 503405-78051.}}}

\author{Liangzu Peng\thanks{Mathematical Institute for Data Science, Johns Hopkins University
  (\email{lpeng25,mahyarfazlyab,rvidal@jhu.edu}).}
\and Mahyar Fazlyab\footnotemark[2]
\and Ren\'e Vidal\footnotemark[2]}

\usepackage{amsopn}

\usepackage{tikz}
\usepackage{comment}
\usepackage{bm}
\usepackage{amsmath,amssymb,rotating, mathrsfs}
\usepackage{mathtools}

\usepackage{color}

\usepackage{epsfig}
\usepackage{graphicx}

\usepackage{subfig}

\usepackage{booktabs}

\usepackage{cancel} 

\usepackage{xfrac} 
\usepackage{xparse} 

\DeclareMathAlphabet{\mathpzc}{OT1}{pzc}{m}{it}
\usepackage{enumitem}

\newcommand{\cI}{\mathcal{I}}

\newcommand{\cN}{\mathcal{N}}

\newcommand{\pA}{\mathpzc{A}}
\newcommand{\pB}{\mathpzc{B}}

\newcommand{\pD}{\mathpzc{D}}

\newcommand{\pL}{\mathpzc{L}}

\newcommand{\pQ}{\mathpzc{Q}}

\newcommand{\pW}{\mathpzc{W}}

\newcommand{\pY}{\mathpzc{Y}}

\newcommand{\pa}{\mathpzc{a}}

\newcommand{\pw}{\mathpzc{w}}

\newcommand{\pz}{\mathpzc{z}}

\newcommand{\bbR}{\mathbb{R}}
\newcommand{\bbS}{\mathbb{S}}

\newcommand{\bC}{\bm{C}}

\newcommand{\bE}{\bm{E}}

\newcommand{\bI}{\bm{I}}

\newcommand{\bP}{\bm{P}}
\newcommand{\bQ}{\bm{Q}}
\newcommand{\bR}{\bm{R}}
\newcommand{\bS}{\bm{S}}
\newcommand{\bT}{\bm{T}}

\newcommand{\bV}{\bm{V}}

\newcommand{\bX}{\bm{X}}

\newcommand{\ba}{\bm{a}}

\newcommand{\bq}{\bm{q}}
\newcommand{\bs}{\bm{s}}

\newcommand{\bw}{\bm{w}}
\newcommand{\bx}{\bm{x}}
\newcommand{\by}{\bm{y}}
\newcommand{\bz}{\bm{z}}

\newcommand{\beps}{\bm{\epsilon}}

\newcommand{\balpha}{\bm{\alpha}}

\NewDocumentCommand{\norm}{mG{2}}{\big\|#1\big\|_{#2}}

\newcommand{\trans}{\top}
\newcommand{\trsp}[1]{#1^\trans}

\DeclareMathOperator{\SO}{SO}

\NewDocumentCommand{\seqp}{mG{n}}{{#1}_1-\cdots+ {#1}_{#2}}
\NewDocumentCommand{\seqm}{mG{n}}{{#1}_1-\cdots- {#1}_{#2}}

\newcommand{\myparagraph}[1]{\smallskip\noindent\textbf{#1.}}

\DeclareMathOperator{\trace}{tr}

\newtheorem{problem}{Problem}

\ifpdf
\hypersetup{
  pdftitle={An Example Article},
  pdfauthor={D. Doe, P. T. Frank, and J. E. Smith}
}
\fi

\begin{document}

\maketitle

\begin{abstract}
The rotation search problem aims to find a 3D rotation that best aligns a given number of point pairs. To induce robustness against outliers for rotation search, prior work considers truncated least-squares (TLS), which is a non-convex optimization problem, and its semidefinite relaxation (SDR) as a tractable alternative. Whether this SDR is theoretically tight in the presence of noise, outliers, or both has remained largely unexplored. We derive conditions that characterize the tightness of this SDR, showing that the tightness depends on the noise level, the truncation parameters of TLS, and the outlier distribution (random or clustered). In particular, we give a short proof for the tightness in the noiseless and outlier-free case, as opposed to the lengthy analysis of prior work.
\end{abstract}



\section{Introduction}\label{section:intro}
	Robust geometric estimation problems in computer vision have been studied for decades \cite{Hartley-2004,Ma-book2004}. However, the analysis of their computational complexity is not sufficiently well understood \cite{Chin-IJCV2020}: There are fast  algorithms that run in real time \cite{Ding-CVPR2020,Bustos-TPAMI2018,Yang-RA-L2020,Peng-CVPR2022}, and there are computational complexity theorems that negate the existence of efficient algorithms \cite{Chin-IJCV2020,Antonante-TRO2021}.\footnote{The catch is that the fast methods might not always be correct (e.g., at extreme outlier rates).} For example, the commonly used \textit{consensus maximization} formulation (for \textit{robust fitting}) is shown to be NP hard in general \cite{Chin-IJCV2020}, and its closely related \textit{truncated least-squares} formulation is not \textit{approximable} \cite{Antonante-TRO2021}, even though they are both highly robust to noise and outliers. Between these ``optimistic'' algorithms and ``pessimistic'' theorems, semidefinite relaxations of truncated least-squares \cite{Lajoie-RAL2019,Yang-ICCV2019,Yang-arXiv2021b} strike a favorable balance between efficiency (as they are typically solvable in polynomial time) and robustness (which is inherited to some extent from the original formulation).  

Even though noise and outliers are ubiquitous in geometric vision, and non-convex formulations and their semidefinite relaxations have been widely used in a large body of papers \cite{Kahl-IJCV2007,Aholt-ECCV2012,Fredriksson-ACCV2012,Cheng-CVPR2015,Maron-ToG2016,Carlone-IROS2015,Carlone-ToR2016,Khoo-TIP2016,Briales-IROS2016,Briales-CVPR2017,Briales-CVPR2018,Probst2019-CVPR2019,Giamou-RAL2019,Agostinho-arXiv2019v2,Li-ICIRS2020,Zhao-CVPR2020,Garcia-IVC2021,Garcia-JMIV2021,Alfassi-Sensors2021,Shi-RSS2021}, much fewer works \cite{Chaudhury-SIAM-J-O2015,Ozyesil-SIAM-J-IS2015,Eriksson-CVPR2018,Rosen-IJRR2019,Iglesias-CVPR2020,Wise-MFI2020,Zhao-TPAMI2020,Tian-TRO2021,Ling-arXiv2021}\footnote{ \cite{Candes-CPAM2013,Bandeira-MP2017,Zhong-SIAM-J-O2018,Lu-SIAM-J-O2019} analyzed SDRs under noise but they are not for geometric vision problems.} provide theoretical insights on the robustness of semidefinite relaxations to noise, only a few semidefinite relaxations \cite{Carlone-RAL2018,Lajoie-RAL2019,Yang-ICCV2019,Yang-arXiv2021b} are empirically robust to outliers, and only one paper \cite{Wang-IMA2013} on \textit{rotation synchronization} gives theoretical guarantees for noise, outliers, and both. Complementary to the story of \cite{Chin-IJCV2020} and inheriting the spirit of \cite{Wang-IMA2013}, in this paper we question \textit{whether ``a specific semidefinite relaxation'' for ``robust rotation search'' is ``tight'' or not}, and provide tightness characterizations that account for the presence of noise, outliers, and both. 

More formally, in this paper we consider the following problem (see \cite{Bustos-PAMI16,Bustos-TPAMI2018,Yang-ICCV2019} for what has motivated this problem):
	\begin{problem}[\textit{Robust Rotation Search}]\label{problem:RRS}
	Let $\{(\by_i, \bx_i)\}_{i=1}^\ell$ be a collection of $\ell$ $3$D point pairs. Assume that a subset $\cI^* \subseteq \{1,\dots,\ell\}$ of these pairs are related by a $3$D rotation $\bR^*_0\in\SO(3)$ up to bounded noise $\{\beps_i:\norm{\beps_i}\leq \delta\}_{i=1}^\ell\subset \bbR^3$ with $\delta\geq 0$, i.e., 
	\begin{align}\label{eq:RRS-model}
		\begin{cases}
			\by_i = \bR^*_0\bx_i + \beps_i, \ \ \ \   & \text{ $i\in \cI^*$ } \\
			\text{ $\by_i$ and $\bx_i$ are arbitrary } & \text{ $i\notin \cI^*$ }.
		\end{cases}
	\end{align}
	Here, $\cI^*$ is called the \textit{inlier} index set. If $i\in \cI^*$ then $\bx_i$, $\by_i$, or $(\by_i,\bx_i)$ is called an inlier, otherwise it is called an \textit{outlier}. The goal is to find $\bR^*_0$ and $\cI^*$ from $\{(\by_i, \bx_i)\}_{i=1}^\ell$.
\end{problem}
To solve this problem, we consider the \textit{truncated least-squares} formulation (rotation version), where the hyper-parameter $c_i^2\geq 0$ is called the \textit{truncation parameter}:
\begin{equation}
	\min_{\bR_0\in\SO(3)} \sum_{i=1}^{\ell} \min \Big\{ \norm{\by_i - \bR_0 \bx_i}^2, \ c_i^2 \Big\} \label{eq:TLS-R} \tag{TLS-R}
\end{equation}
While \eqref{eq:TLS-R} is highly robust to outliers and noise \cite{Yang-arXiv2021b}, it is non-convex and hard to solve. Via a remarkable sequence of algebraic manipulations, \cite{Yang-ICCV2019} showed that \eqref{eq:TLS-R} is equivalent to some non-convex \textit{quadratically constrained quadratic program} \eqref{eq:P-QCQP}, which can be relaxed as a \textit{semidefinite program} \eqref{eq:SDP} via the standard  \textit{lifting} technique. The exact forms of \eqref{eq:P-QCQP} and \eqref{eq:SDP} will be shown in \cref{section:review-SDP}. One approach to study how much robustness \eqref{eq:SDP} inherits from \eqref{eq:TLS-R} or \eqref{eq:P-QCQP} is to verify if the solution of \eqref{eq:SDP} leads to a global minimizer to \eqref{eq:P-QCQP}. Informally, if this is true, then we say that \eqref{eq:SDP} is \textit{tight} (cf. Definition \eqref{def:tightness}). Here, we make the following contributions:
\begin{itemize}
	\item For noiseless point sets without outliers ($\beps_i=0,\cI^*=\{1,\dots,\ell\}$ in \cref{problem:RRS}), we prove that \eqref{eq:SDP} is always tight (\cref{theorem:noiseless+outlier-free}). While this result had already been proven in \cite[Section E.3]{Yang-ICCV2019}, our proof is simpler and shorter.
	\item For noiseless point sets with outliers, \cref{theorem:noiseless+outlier+c} states that \eqref{eq:SDP} is tight for sufficiently small truncation parameters $c_i^2$ and \textit{random} outliers (regardless of the number of outliers), but it is not tight if $c_i^2$ is set too large. \cref{theorem:noiseless+outlier+c2} reveals that \eqref{eq:SDP} is vulnerable to (e.g., not tight in the presence of) \textit{clustered} outlier point pairs that are defined by a rotation different from $\bR_0^*$. Different from \cref{theorem:noiseless+outlier-free}, outliers and improper choices of $c_i^2$ might actually undermine the tightness of \eqref{eq:SDP}. 
	
	\item For noisy point sets without outliers, \cref{theorem:noisy+outlier-free2}, with a technical proof, shows that \eqref{eq:SDP} is tight for sufficiently small noise and for sufficiently large $c_i^2$.
	
	\item The case of noisy data with outliers is the most challenging, but from our analysis of the two previous cases, a tightness characterization for this difficult case follows (\cref{theorem:noisy+outliers}). Thus, we will discuss this case only sparingly.
\end{itemize}

\myparagraph{Paper Organization} In \cref{section:review-SDP} we review the derivations of \eqref{eq:SDP} from \eqref{eq:TLS-R} \cite{Yang-ICCV2019}, while we also provide new insights. In \cref{section:main_results}, we discuss our main results. In \cref{section:discussion}, we present limitations of our work and potential avenues for future research.

\myparagraph{Notations and Basics} We employ the MATLAB notation $[a_1;\dots;a_\ell]$ to denote concatenation into a column vector. Given a $4(\ell+1)\times 4(\ell+1)$ matrix $\pA$, we employ the bracket notation $[ \pA ]_{ij}$ of \cite{Yang-ICCV2019} to denote the $4\times 4$ submatrix of $\pA$ whose rows are indexed by $\{4i+1,\dots,4i+4\}$ and columns by $\{4j+1,\dots,4j+4\}$. Following our previous work on robust rotation search \cite{Peng-CVPR2022}, we treat \textit{unit quaternions} as unit vectors on the $3$-sphere $\bbS^3$. It is known that every 3D rotation $\bR\in\SO(3)$ can be equivalently written as
\begin{equation}\label{eq:rotation-quaternion}
	\bR=\begin{bmatrix}
		w_1^2+w_2^2-w_3^2 -w_4^2 & 2(w_2w_3-w_1w_4)& 2(w_2w_4+w_1w_3) \\
		2(w_2w_3+w_1w_4)& w_1^2+w_3^2-w_2^2-w_4^2& 2(w_3w_4-w_1w_2)\\
		2(w_2w_4-w_1w_3)& 2(w_3w_4+w_1w_2)& w_1^2+w_4^2-w_2^2-w_3^2\\
	\end{bmatrix}
\end{equation}
where $\bw=[w_1;w_2;w_3;w_4]\in\bbS^3$ and $-\bw$ are unit quaternions. Conversely, every $3\times 3$ matrix of the form \eqref{eq:rotation-quaternion} with $[w_1;w_2;w_3;w_4]\in\bbS^3$ is a $3$D rotation. This means a two-to-one correspondence between unit quaternions ($\bbS^3$) and $3$D rotations ($\SO(3)$).

\section{\eqref{eq:TLS-R} and Its Relaxation: Review and New Insights}\label{section:review-SDP}
In \cref{subsection:derivation} we derive a semidefinite relaxation \eqref{eq:SDP} from \eqref{eq:TLS-R}. More specifically, we show that \eqref{eq:TLS-R} is equivalent to a truncated least-squares problem named \eqref{eq:TLS-Q}, with the optimization variable being a unit quaternion. We further show that \eqref{eq:TLS-Q} can be equivalently written as a quadratically constrained quadratic program, labeled as \eqref{eq:P-QCQP}. Thus, we obtain the semidefinite relaxation \eqref{eq:SDP} of \eqref{eq:P-QCQP}, as a result of \textit{lifting}, and their dual program \eqref{eq:Dual}, via standard Lagriangian calculation. See Table \ref{table:programs} for an overview.

\setlength{\tabcolsep}{0.3pt}
\renewcommand{\arraystretch}{1.5}
\begin{table}[h]
	\centering
	\caption{Descriptions of different programs. TLS means Truncated Least-Squares. \label{table:programs}}
	\begin{tabular}{cc}
		\toprule
		Programs & Description \\
		\midrule
		\eqref{eq:TLS-R} & The TLS problem with the optimization variable being a 3D rotation\\ 
		\eqref{eq:TLS-Q} & \ \ TLS with the optimization variable being a unit quaternion, equivalent to \eqref{eq:TLS-R} \\
		\eqref{eq:P-QCQP} & A quadratically constrained quadratic program, equivalent to \eqref{eq:TLS-Q} \\
		\midrule
		\eqref{eq:SDP} & The semidefinite relaxation of \eqref{eq:P-QCQP}, obtained via \textit{lifting} \\
		\eqref{eq:Dual} & The dual program of \eqref{eq:SDP} and \eqref{eq:P-QCQP} \\
		\bottomrule
	\end{tabular}
\end{table}
While \cref{subsection:derivation} follows the development of \cite{Yang-ICCV2019} in spirit,  our derivation is simpler. For example, we dispensed with the use of \textit{quaternion product}  \cite{Horn-JOSAA1987} in \cite{Yang-ICCV2019}, which is a sophisticated algebraic operation. That said, it is safe to treat our \eqref{eq:SDP} as equivalent to the \textit{naive relaxation} of \cite{Yang-ICCV2019}; see the appendix for a detailed discussion. 

In \cref{subsection:optimality-conditions} we discuss the KKT optimality conditions that are essential for studying the interplay among \eqref{eq:P-QCQP}, \eqref{eq:SDP}, and \eqref{eq:Dual} and the tightness of \eqref{eq:SDP}.

\subsection{Derivation}\label{subsection:derivation} 
The term $\| \by_i - \bR_0 \bx_i \|_2^2=\| \by_i \|_2^2 + \| \bx_i\|_2^2 - 2\trsp{\by_i} \bR_0\bx_i$ of \eqref{eq:TLS-R} depends linearly on the rotation $\bR_0$. Moreover, each entry of any rotation $\bR_0$ depends quadratically on its \textit{unit quaternion} representation $\bw_0\in\bbS^3$; recall \eqref{eq:rotation-quaternion}. One then naturally asks whether $\| \by_i - \bR_0 \bx_i \|_2^2$ is a quadratic form in $\bw_0$; the answer is affirmative:
\begin{lemma}[Rotations and Unit Quaternions]\label{lemma:rotation-quaternion}
	Let $\bR_0$ be a 3D rotation, then we have
	\begin{equation}
		\norm{\by_i-\bR_0\bx_i}^2 = \trsp{\bw_0} \bQ_i \bw_0,
	\end{equation}
	where $\bQ_i$ is a $4\times 4$ positive semidefinite matrix, and $\bw_0\in\bbS^3$ is the unit quaternion representation of $\bR_0$. Moreover, the eigenvalues of $\bQ_i$ are respectively
	\begin{equation}\label{eq:eigenvalues-rotation-quaternion}
		\big(\norm{\by_i}+\norm{\bx_i}\big)^2, \big(\norm{\by_i}+\norm{\bx_i}\big)^2, \big(\norm{\by_i}-\norm{\bx_i}\big)^2, \big(\norm{\by_i}-\norm{\bx_i}\big)^2.
	\end{equation}
\end{lemma}
\begin{proof}[Proof of \cref{lemma:rotation-quaternion}]
	To simplify notations, here, and only here, we drop all indices and consider any 3D rotation $\bR$ and any 3D point pairs $(\by,\bx)$. We will first show
	\begin{equation}\label{eq:Q(y,x)}
		\norm{\by-\bR\bx}^2 = \trsp{\bw}\ Q(\by,\bx) \ \bw,
	\end{equation}
	where $\bw$ is the unit quaternion representation of $\bR$, and $Q:\bbR^3\times \bbR^3\to \bbR^{4\times 4}$ is some function that maps $(\by,\bx)$ to a $4\times 4$ matrix $Q(\by,\bx)$. Recall \eqref{eq:rotation-quaternion}, and we note that the first, second, and third entries of $\bR \bx$ can be written as $\trsp{\bw}\bX_1\bw$, $\trsp{\bw}\bX_2\bw$, and $\trsp{\bw}\bX_3\bw$, respectively, where $\bX_1$, $\bX_2$, $\bX_3$ are symmetric  matrices respectively defined as:
	\begin{equation*}
		\small{ \begin{matrix}
				\bX_1 & \bX_2 & \bX_3 \\
				\begin{bmatrix}
					x_1&0&x_3&-x_2\\
					0&x_1&x_2&x_3\\
					x_3&x_2&-x_1&0\\
					-x_2&x_3&0& -x_1
				\end{bmatrix} &
				\begin{bmatrix}
					x_2&-x_3&0&x_1\\
					-x_3&-x_2&x_1&0\\
					0&x_1&x_2&x_3\\
					x_1&0&x_3&-x_2
				\end{bmatrix}&
				\begin{bmatrix}
					x_3&x_2&-x_1&0\\
					x_2&-x_3&0&x_1\\
					-x_1&0&-x_3&x_2\\
					0&x_1&x_2&x_3
				\end{bmatrix} 
			\end{matrix}
		}
	\end{equation*}
	Define $U(\by,\bx):=y_1\bX_1+y_2\bX_2 + y_3\bX_3$, then $\trsp{\by}\bR\bx=\trsp{\bw}\ U(\by,\bx)\ \bw$. With
	\begin{equation}\label{eq:def-Qi}
		Q(\by,\bx) := \big(\norm{\by}^2 + \norm{\bx}^2\big)\bI_4-2U(\by,\bx), 
	\end{equation}
	we have proved \eqref{eq:Q(y,x)}. Note that, as the above derivation implies, given any unit quaternion $\bw'\in\bbS^3$, there is a unique 3D rotation $\bR'$ which satisfies $\trsp{(\bw')} \ Q(\by,\bx)\  \bw'=\| \by-\bR'\bx \|_2^2$; the latter term $\norm{\by-\bR'\bx}^2$ is always nonnegative. Combining this with the fact that $Q(\by,\bx)$ is symmetric, we get that $Q(\by,\bx)$ is positive semidefinite.
	
	It remains to find the eigenvalues of $Q(\by,\bx)$. Note that to minimize  $\trsp{\bw}\ Q(\by,\bx) \ \bw$ over $\bw\in\bbS^3$ is to minimize $\| \by - \bR \bx \|_2^2$ over $\bR\in\SO(3)$. Combining this with a geometric fact gives
	\begin{equation}
		\min_{\bw\in\bbS^3} \trsp{\bw}\ Q(\by,\bx) \ \bw = \min_{\bR\in\SO(3)} \norm{\by - \bR \bx}^2 = \big(\norm{\by}-\norm{\bx}\big)^2, \label{eq:min_residual} 
	\end{equation}
	So $Q(\by,\bx)$ has $(\| \by \|_2 - \| \bx \|_2)^2$ as its minimum eigenvalue. Since the minimum $(\| \by \|_2 - \| \bx \|_2)^2$ of \eqref{eq:min_residual} is attained at (at least) two different rotations $\bR$'s which make $\bR\bx$ and $\by$ parallel and pointing to the same direction, i.e., there are (at least) two different unit quaternions $\pm\bw$ and $\pm\bw'$ corresponding to the minimum eigenvalue $(\| \by \|_2 - \| \bx \|_2)^2$. Thus the eigenspace of $Q(\by,\bx)$ associated with eigenvalue $(\| \by \|_2 - \| \bx \|_2)^2$ is of dimension at least $2$. Repeat the above argument for the maximum eigenvalue $(\| \by \|_2 + \| \bx \|_2)^2$, and note that $Q(\by,\bx)$ is of size $4\times 4$, then we know all eigenvalues of $Q(\by,\bx)$. With $\bQ_i:=Q(\by_i,\bx_i)$, we finish the proof.
\end{proof}

While the exact form of $\bQ_i$ is complicated, \cref{lemma:rotation-quaternion} informs us of a characterization on eigenvalues of $\bQ_i$, which is much easier to work with. Note that, while the relation between 3D rotations and unit quaternions is well known (see, e.g., \cite{Horn-JOSAA1987}), we have not found \eqref{eq:eigenvalues-rotation-quaternion} in the literature, except in the appendix of our prior work \cite{Peng-CVPR2022}.

From \cref{lemma:rotation-quaternion}, we now see that \eqref{eq:TLS-R} is equivalent to
\begin{equation}
	\min_{\bw_0\in \bbS^3} \sum_{i=1}^{\ell} \min \Big\{ \trsp{\bw}_0\bQ_i \bw_0,\ c_i^2 \Big\}. \label{eq:TLS-Q} \tag{TLS-Q}
\end{equation}
Using the following simple equality ($\theta\in\{-1,1\}$ in \cite{Lajoie-RAL2019,Yang-ICCV2019}; see also \cite{Ikami-CVPR2018})
\begin{equation*}
	\min\{ a,b \}= \min_{\theta\in\{0, 1\}} \theta a + (1-\theta) b= \min_{\theta^2=\theta} \theta a + (1-\theta) b,
\end{equation*}
problem \eqref{eq:TLS-Q} can be equivalently written as
\begin{equation}
	\min_{\bw_0\in\bbS^3,\theta_i^2=\theta_i} \sum_{i=1}^{\ell}\Big( \theta_i\trsp{\bw_0}\bQ_i \bw_0 - \theta_i c_i^2 \Big) + \sum_{i=1}^{\ell} c_i^2 \label{eq:TLS-wtheta}
\end{equation}
Note that, while the constant $\sum_{i=1}^{\ell} c_i^2$ in \eqref{eq:TLS-wtheta} can be ignored, keeping it there will simplify matters. Even though the objective of \eqref{eq:TLS-wtheta} is a cubic polynomial in entries of the unknowns $\bw_0$ and $\theta_i$'s, problem \eqref{eq:TLS-wtheta} is equivalent to a quadratic program. Indeed, let $\bw_i:=\theta_i \bw_0$, which implies  $\theta_i=\trsp{\bw_0}\bw_i$. Then \eqref{eq:TLS-wtheta} becomes
\begin{equation}\label{eq:QCQP}
	\min_{\bw_0\in\bbS^3,\bw_i\in \{\bw_0, 0\}}  \sum_{i=1}^{\ell}\Big( \trsp{\bw_0}(\bQ_i - c_i^2\bI_4) \bw_i \Big) + \sum_{i=1}^{\ell} c_i^2
\end{equation}
Problem \eqref{eq:QCQP} now has its objective function quadratic, while in fact its constraints are also quadratic. To see this, one easily verifies that the binary constraint $\bw_i\in\{\bw_0, 0\}$ can be equivalently written quadratically as $\bw_i\trsp{\bw_0} = \bw_i\trsp{\bw_i}$. Collecting all vectors of variables into a $4(\ell+1)$ dimensional column vector $\pw=[\bw_0;\dots;\bw_\ell]$, we have
\begin{equation}\label{eq:equivalent-constraints}
	\begin{cases}
		\bw_0\in\bbS^3 \\
		\bw_i\in \{\bw_0, 0\}
	\end{cases} \Leftrightarrow
	\begin{cases}
		\trace(\bw_0\trsp{\bw_0})=1 \\
		\bw_i\trsp{\bw_0} = \bw_i\trsp{\bw_i}
	\end{cases}\Leftrightarrow
	\begin{cases}
		\trace\big([\pw\trsp{\pw}]_{00}\big) =  1 \\
		[\pw\trsp{\pw}]_{0i} = [\pw\trsp{\pw}]_{ii}
	\end{cases}
\end{equation}
In the last equivalence of \eqref{eq:equivalent-constraints} we used the notation $[\ \cdot\ ]_{ij}$ of \cref{section:intro}. Having confirmed \eqref{eq:equivalent-constraints}, we can now equivalently transform \eqref{eq:QCQP} into the following \eqref{eq:P-QCQP}:
\begin{equation}\label{eq:P-QCQP}
	\begin{split}
		\min_{\pw\in\bbR^{4(\ell+1)}} &\ \ \trace\Big(\pQ \pw \trsp{\pw}\Big) + \sum_{i=1}^{\ell} c_i^2  \\
		\text{s.t.}& \ \ [\pw\trsp{\pw}]_{0i} = [\pw\trsp{\pw}]_{ii}, \ \ \forall\ i\in\{1,\dots,\ell\} \\
		& \ \ \trace\big([\pw\trsp{\pw}]_{00}\big) =  1 
	\end{split} \tag{QCQP}
\end{equation}
In \eqref{eq:P-QCQP}, $\pQ$ is our $4(\ell+1)\times 4(\ell+1)$ data matrix, symmetric and satisfying
\begin{equation}\label{eq:def-tildepQ}
	\begin{cases}
		[\pQ]_{0i}=[\pQ]_{i0}=\frac{1}{2}(\bQ_{i}-c_i^2\bI_4),\ \ \forall\ i\in\{1,\dots,\ell\} \\
		\text{all other entries of $\pQ$ are zero.}
	\end{cases} 
\end{equation}

It is now not hard to derive the semidefinite relaxation \eqref{eq:SDP} and dual program \eqref{eq:Dual} from \eqref{eq:P-QCQP} via lifting and standard Lagrangian calculation respectively:

\begin{lemma}[\eqref{eq:SDP} and \eqref{eq:Dual}]\label{lemma:SDP-D}
	The dual and semidefinite relaxation of \eqref{eq:P-QCQP} are
	\begin{equation}\label{eq:Dual}
		\max_{\mu,\pD}\ \ \mu  + \sum_{i=1}^{\ell} c_i^2 \ \ \ \ \ \text{s.t.} \ \ \ \ \ \pQ-\mu\pB-\pD \succeq 0 \tag{D}
	\end{equation}
	\begin{equation}
		\begin{split}
			\min_{\pW\succeq 0} &\ \ \trace\big(\pQ \pW\big)  + \sum_{i=1}^{\ell} c_i^2  \\
			\text{s.t.}& \ \ [\pW]_{0i} = [\pW]_{ii},\ \ \forall\ i\in\{1,\dots,\ell\},  \ \ \ \ \trace\big([\pW]_{00}\big) = 1 
		\end{split} \tag{SDR} \label{eq:SDP}
	\end{equation}
	In \eqref{eq:Dual}, $\pB\in\bbR^{4(\ell+1)\times 4(\ell+1)}$ is a matrix of zeros except $[\pB]_{00}=\bI_4$, while $\pD\in\bbR^{4(\ell+1)\times 4(\ell+1)}$ is a matrix of dual variables accounting for the $\ell$ constraints $[\pw\trsp{\pw}]_{0i} = [\pw\trsp{\pw}]_{ii}$, i.e., $\pD$ satisfies 
	\begin{equation}\label{eq:D-constraint}
		\begin{cases}
			\text{$\pD$ is symmetric}, \ \ [\pD]_{ii} +2[\pD]_{0i}=0, \ \ \forall\ i\in\{1,\dots,\ell\} \\
			\text{all other entries of $\pD$ are zero.}
		\end{cases}
	\end{equation}
\end{lemma}
\begin{proof}[Proof of \cref{lemma:SDP-D}]
	\cref{lemma:SDP-D} is fairly standard, and we provide a concise proof for completeness. The Lagrangian of \eqref{eq:P-QCQP} is given as
	\begin{equation*}
		\pL_{\textnormal{\ref{eq:P-QCQP}}}(\pw,\mu, \pD) = \trace\Big((\pQ-\mu \pB-\pD) \pw \trsp{\pw}\Big) + \mu + \sum_{i=1}^{\ell}c_i^2,
	\end{equation*}
	where $\mu\in\bbR$ is the Lagrangian multiplier that accounts for constraint $\trace\big([\pw\trsp{\pw}]_{00}\big) =  1$, and $\pB$ and $\pD\in\bbR^{4(\ell+1)\times 4(\ell+1)}$ have the said form. The dual function $g_{\text{D}}(\mu, \pD)$ then is
	\begin{equation*}
		g_{\textnormal{\ref{eq:Dual}}}(\mu,\pD) = \inf_{\pw} \pL_{\textnormal{\ref{eq:P-QCQP}}}(\pw,\mu, \pD)=\begin{cases}
			\mu + \sum_{i=1}^{\ell}c_i^2 & \pQ-\mu \pB - \pD \succeq 0\\
			-\infty& \text{ otherwise}
		\end{cases}
	\end{equation*}
	The dual problem of \eqref{eq:P-QCQP} is thus  \eqref{eq:Dual}, where we maximize $g_\textnormal{\ref{eq:Dual}}$ over the dual variables $\mu$ and $\pD$. Next, we show that \eqref{eq:SDP} is the dual of \eqref{eq:Dual}. The Lagrangian of \eqref{eq:Dual} is
	\begin{equation*}
		\pL_{\textnormal{\ref{eq:Dual}}}(\pW,\mu, \pD) = \trace\Big((\pQ-\mu \pB-\pD) \pW\Big) + \mu + \sum_{i=1}^{\ell}c_i^2.
	\end{equation*}
	For weak duality to hold, that is for $\pL_{\textnormal{\ref{eq:Dual}}}(\pW,\mu, \pD)$ to majorize $g_{\textnormal{\ref{eq:Dual}}}(\mu,\pD)$ at any feasible dual points $(\mu, \pD)$ with $\pQ-\mu \pB - \pD \succeq 0$, we need to add an extra constraint $\pW\succeq 0$. Then, maximizing $\pL_{\textnormal{\ref{eq:Dual}}}(\pW,\mu, \pD)$ for any $\mu\in\bbR$ and any $\pD$ of the form \eqref{eq:D-constraint} and then minimizing the resulting dual function of \eqref{eq:Dual} over $\pW\succeq 0$ gives \eqref{eq:SDP} as the dual program of \eqref{eq:Dual}. 
\end{proof}

\subsection{The Tightness and Optimality Conditions}\label{subsection:optimality-conditions}
We start with the following definition:
\begin{definition}[Tightness]\label{def:tightness}
	\eqref{eq:SDP} is said to be \textit{tight} if it admits $\hat{\pw} \trsp{(\hat{\pw})}$ as a global minimizer, where $\hat{\pw}\in\bbR^{4(\ell+1)}$ globally minimizes \eqref{eq:P-QCQP}.
\end{definition}
We now proceed with some basic observations. First, weak duality between \eqref{eq:P-QCQP} and \eqref{eq:Dual} holds as a result of Lagrangian calculation. Second:
\begin{proposition}\label{prop:strong-duality}
	Strong duality between \eqref{eq:SDP} and \eqref{eq:Dual} holds.
\end{proposition}
\begin{proof}
	It suffices to exhibit a \textit{strictly feasible} point $(\mu,\pD)$ of \eqref{eq:Dual} (such that $\pQ-\mu\pB-\pD \succ 0$) and then invoke Slater's condition \cite{Boyd-book2004}. Let $\pD$ be of the form \eqref{eq:D-constraint} satisfying $[\pD]_{0i}:=\frac{1}{2}(\bQ_i +c_i^2\bI_4)$. Then $\pQ-\mu\pB-\pD \succ 0$ is equivalent to  
	\begin{equation*}
		\begin{split}
			&\ \ \trsp{[\bz_0;\dots;\bz_\ell]} \big( \pQ-\mu \pB-\pD \big) [\bz_0;\dots;\bz_\ell] > 0, \ \ \forall \bz_i\in\bbR^4 \\
			\Leftrightarrow&\ \ -\mu \norm{\bz_0}^2 -  \sum_{i=1}^{\ell} c_i^2\norm{\bz_0}^2 +  \sum_{i=1}^{\ell} \Big(c_i^2\norm{\bz_0 - \bz_i}^2 + \trsp{\bz_{i}}  \bQ_i \bz_i \Big) > 0,\ \  \forall \bz_i\in\bbR^4 
		\end{split}
	\end{equation*}
	Since  \cref{lemma:rotation-quaternion} implies $\bQ_i\succeq 0$, the above holds for sufficiently small $\mu$. Such $(\mu,\pD)$ is strictly feasible, thus finishing the proof.
\end{proof}
Since the strong duality between \eqref{eq:SDP} and \eqref{eq:Dual} holds and every feasible point of \eqref{eq:QCQP} is also feasible for \eqref{eq:SDP}, we can then obtain the following necessary and sufficient optimality condition for characterizing the tightness of \eqref{eq:SDP}:

\begin{proposition}[Necessary and Sufficient Optimality Conditions]\label{prop:tightnesscondition}
	Suppose that \eqref{eq:TLS-Q} preserves all inliers $\bQ_1,\dots,\bQ_{k^*}$ and rejects all outliers $\bQ_{k^*+1},\dots,\bQ_{\ell}$ outliers. Recall that $\hat{\bw}_0$ denotes a global minimizer of \eqref{eq:TLS-Q}. \eqref{eq:SDP} is tight if and only if there is a matrix $\hat{\pD}$ of the form \eqref{eq:D-constraint} that satisfies the following conditions:
	\begin{equation}
		\hspace{-0.65cm} \textnormal{Stationarity Condition:}  \begin{cases}
			\big(2[\hat{\pD}]_{0i} + \bQ_i-c_i^2\bI_4\big)\hat{\bw}_0=0,  & \ \forall i\in\{1,\dots,k^* \} \\
			\big(2[\hat{\pD}]_{0j}+c_j^2\bI_4 - \bQ_j\big)\hat{\bw}_0=0,  &\ \forall j\in \{k^*+1,\dots,\ell \} \tag{O1}
		\end{cases} \label{eq:stationary2} \\
	\end{equation}
	\begin{equation}
		\textnormal{Dual Feasibility: } -\hat{\mu} \norm{\bz_0}^2 + 2\sum_{i=1}^{\ell}\trsp{\bz_i}[\hat{\pD}]_{0i} \bz_i -\sum_{i=1}^{\ell}\trsp{\bz_0}\big(2[\hat{\pD}]_{0i} - \bQ_i+c_i^2\bI_4\big)\bz_i\geq 0, \forall \bz_i\in\bbR^4 \tag{O2} \label{eq:dualfeasibility2}
	\end{equation}
	\begin{equation}
		\hspace{-4.7cm} \textnormal{Objective Value Exactness:\ \ } \hat{\mu}= \sum_{i=1}^{k^*} \big(\trsp{\hat{\bw}_0}\bQ_i \hat{\bw}_0 - c_i^2\big) \tag{O3} \label{eq:value-mu2}
	\end{equation}
\end{proposition}

\begin{proof}[Proof of \cref{prop:tightnesscondition}]
	Since \eqref{eq:TLS-Q} preserves all inliers $\bQ_1,\dots,\bQ_{k^*}$ and rejects all outliers, $\hat{\pw}=[\hat{\bw}_0;\dots; \hat{\bw}_0; 0;\dots,0]$, where $\hat{\bw}_0$ appeared $k^*+1$ times, is a global minimizer of \eqref{eq:P-QCQP} and reaches the minimum value $\sum_{i=1}^{k^*} \trsp{\hat{\bw}_0}\bQ_i \hat{\bw}_0  + \sum_{j=k^*+1}^{\ell} c_j^2$. Clearly, $\hat{\pw} \trsp{\hat{\pw}}$ is a feasible point of \eqref{eq:SDP}. Since strong duality between \eqref{eq:SDP} and \eqref{eq:Dual} holds (\cref{prop:strong-duality}), the minimum of \eqref{eq:SDP} is $\hat{\mu} + \sum_{i=1}^\ell c_i^2$, and \eqref{eq:SDP} is tight (i.e., $\hat{\pw} \trsp{\hat{\pw}}$ globally minimizes \eqref{eq:SDP}) if and only if (the minima of \eqref{eq:SDP} and \eqref{eq:P-QCQP} coincide) the KKT conditions
	\begin{equation*}
		(\pQ-\hat{\mu} \pB-\hat{\pD})\hat{\pw} \trsp{\hat{\pw}} = 0, \ \ \ \ \  
		\pQ-\hat{\mu} \pB-\hat{\pD}\succeq 0
	\end{equation*}
	hold true. via basic algebra and calculation, we obtain that \eqref{eq:SDP} is tight if and only if
	\begin{equation*}
		(\pQ-\hat{\mu} \pB-\hat{\pD})\hat{\pw}  = 0, \ \ \ \ \  
		\pQ-\hat{\mu} \pB-\hat{\pD}\succeq 0,\ \ \  \hat{\mu}= \sum_{i=1}^{k^*} \big(\trsp{\hat{\bw}_0}\bQ_i \hat{\bw}_0 - c_i^2\big).
	\end{equation*}
	With $\hat{\mu}= \sum_{i=1}^{k^*} \big(\trsp{\hat{\bw}_0}\bQ_i \hat{\bw}_0 - c_i^2\big)$, the first condition $(\pQ-\hat{\mu} \pB-\hat{\pD})\hat{\pw}  = 0$ is equivalent to 
	\begin{equation}\label{eq:stationary-in-proof}
		\begin{cases}
			\sum_{i=1}^{k^*} \frac{\bQ_i-c_i^2\bI_4}{2} \hat{\bw}_0 - \sum_{i=1}^{k^*} \big(\trsp{\hat{\bw}_0}\bQ_i \hat{\bw}_0 - c_i^2\big)\hat{\bw}_0  - \sum_{i=1}^{k^*} [\pD]_{0i} \hat{\bw}_0 = 0 \\
			\frac{\bQ_i-c_i^2\bI_4}{2} \hat{\bw}_0 - [\pD]_{0i} \hat{\bw}_0 - [\pD]_{ii} \hat{\bw}_0 = 0, \ \ \forall i\in\{ 1,\dots,k^* \}  \\
			\frac{\bQ_j-c_j^2\bI_4}{2} \hat{\bw}_0 - [\pD]_{0j} \hat{\bw}_0	= 0, \ \ \forall j\in\{ k^*+1,\dots,\ell \}
		\end{cases}
	\end{equation}
	Note that our assumption on \eqref{eq:TLS-Q} implies that  $\hat{\bw}_0$ is an eigenvector of $\sum_{i=1}^{k^*} \bQ_i$ corresponding to its minimum eigenvalue, i.e., $\sum_{i=1}^{k^*} \bQ_i\hat{\bw}_0 =  \big(\sum_{i=1}^{\ell} \trsp{\hat{\bw}_0} \bQ_i \hat{\bw}_0 \big)  \hat{\bw}_0 $. Also note that $\pD$ is symmetric with $[\pD]_{ii} +2[\pD]_{0i}=0$ \eqref{eq:D-constraint}. The above \eqref{eq:stationary-in-proof} can then be simplified into 
	\begin{equation*}
		\begin{cases}
			\sum_{i=1}^{k^*} \frac{\bQ_i-c_i^2\bI_4}{2} \hat{\bw}_0 + \sum_{i=1}^{k^*} [\pD]_{0i} \hat{\bw}_0 = 0 \\
			\frac{\bQ_i-c_i^2\bI_4}{2} \hat{\bw}_0 + [\pD]_{0i} \hat{\bw}_0  = 0, \ \ \forall i\in\{ 1,\dots,k^* \}  \\
			\frac{\bQ_j-c_j^2\bI_4}{2} \hat{\bw}_0 - [\pD]_{0j} \hat{\bw}_0	= 0, \ \ \forall j\in\{ k^*+1,\dots,\ell \}
		\end{cases}
	\end{equation*}
	which is equivalent to \eqref{eq:stationary2}. The second condition $\pQ-\hat{\mu} \pB-\hat{\pD}\succeq 0$
	is equivalent to 
	\begin{equation*}
		\trsp{[\bz_0;\dots;\bz_\ell]} \big( \pQ-\hat{\mu} \pB-\hat{\pD} \big) [\bz_0;\dots;\bz_\ell] \geq 0, \ \ \forall \bz_i\in\bbR^4,
	\end{equation*}
	which is the same as \eqref{eq:dualfeasibility2} by rewriting using the definitions of $\pQ$ \eqref{eq:def-tildepQ} and $\pD$ \eqref{eq:D-constraint}.
\end{proof}

Armed with \cref{prop:tightnesscondition}, to show \eqref{eq:SDP} is tight or not we need to find the \textit{dual certificates} $[\hat{\pD}]_{0i}$'s (and $\hat{\mu}$) that fulfill the (simplified) optimality conditions. Identifying eligible $[\hat{\pD}]_{0i}$'s or showing that such $[\hat{\pD}]_{0i}$'s do not exist is a core idea in proving our main results, which are to be discussed in greater detail in the next section. 

\section{Tightness Characterizations}\label{section:main_results}
We now present our results regarding the tightness of \eqref{eq:SDP}.  \Cref{subsection:noiseless+outlier-free} treats the simplest noiseless + outlier-free case (\cref{theorem:noiseless+outlier-free}). \Cref{subsection:noiseless+outlier+c,subsection:noisy+outlier-free} take outliers (\cref{theorem:noiseless+outlier+c,theorem:noiseless+outlier+c2}) and noise (\cref{theorem:noisy+outlier-free2}) into account respectively; \cref{subsection:noisy+outliers} brings them together for the noisy + outliers case (\cref{theorem:noisy+outliers}).

\subsection{The Noiseless + Outlier-Free Case}\label{subsection:noiseless+outlier-free}
\begin{theorem}[Noiseless and Outlier-Free Point Sets]\label{theorem:noiseless+outlier-free}
	\eqref{eq:SDP} is tight in the absence of noise and outliers, meaning that $\pw^*\trsp{(\pw^*)}$ globally minimizes \eqref{eq:SDP}, where $\pw^*=[\bw_0^*;\dots;\bw_0^*]\in\bbR^{4(\ell+1)}$ is a global minimizer of \eqref{eq:P-QCQP}.
\end{theorem}
\begin{proof}
	Note that $\pw^*=[\bw^*_0;\dots;\bw^*_0]$ is a global minimizer of \eqref{eq:P-QCQP} that results in the optimal value $0$. Let $\hat{\pD}$ satisfy the constraint \eqref{eq:D-constraint} with $[\hat{\pD}]_{0i}:= \frac{1}{2}(\bQ_i +c_i^2\bI_4)$ for every $i=1,\dots,\ell$ and let $\hat{\mu}:=-\sum_{i=1}^\ell c_i^2$. Then, with $\bQ_i\bw^*_0=0$ (\cref{lemma:rotation-quaternion}), one easily verifies that optimality conditions \eqref{eq:stationary2} and \eqref{eq:value-mu2} of \cref{prop:tightnesscondition} hold. It remains to prove condition \eqref{eq:dualfeasibility2}. Substitute the values of  $[\hat{\pD}]_{0i}$, $\hat{\mu}$ into \eqref{eq:dualfeasibility2} and it simplifies:
	\begin{equation}
		\begin{split}
			&\sum_{i=1}^{\ell} c_i^2\norm{\bz_0}^2+\sum_{i=1}^{\ell}\trsp{\bz_{i}} \big( \bQ_i+c_i^2\bI_4 \big) \bz_i -  2\sum_{i=1}^{\ell} c_i^2\trsp{\bz_i} \bz_0\geq 0, \ \ \forall \bz_i\in\bbR^4 \\
			\Leftrightarrow\ \ & \sum_{i=1}^{\ell} \Big(c_i^2\norm{\bz_0 - \bz_i}^2 + \trsp{\bz_{i}}  \bQ_i \bz_i \Big) \geq 0, \ \ \forall \bz_i\in\bbR^4
		\end{split}
	\end{equation}
	Thus \eqref{eq:dualfeasibility2} holds, as every $\bQ_i$ is positive semidefinite (\cref{lemma:rotation-quaternion}). One also observes that the equality is attained if and only if $\bz_0=\cdots=\bz_\ell=\bw_0^*$ or $\bz_0=\cdots=\bz_\ell=0$.
\end{proof}
Our contribution here is a shorter proof for \cref{theorem:noiseless+outlier-free} than \cite{Yang-ICCV2019}. Besides \cref{lemma:rotation-quaternion}, another key idea that shortens the proof is our construction of the dual certificate $\hat{\pD}$ (or $[\hat{\pD}]_{0i}$'s). While constructing dual certificates might be an art as there might not exist general approaches for doing so, our experience is to (1) start with the \textit{simplest} case (e.g., noiseless + outlier-free), (2) make \textit{observations}: observe the optimality conditions (cf. \cref{prop:tightnesscondition}), inspect the first and second order Riemannian optimality conditions (cf. \cite{Boumal-2020}), discover some properties of data (e.g., \cref{lemma:rotation-quaternion}), (3) \textit{repeatedly try} different choices of certificates. 

\subsection{The Noiseless + Outliers Case}\label{subsection:noiseless+outlier+c}
Different from Theorem \ref{theorem:noiseless+outlier-free}, the tightness of \eqref{eq:SDP} in the presence of outliers depends on the truncation parameters $c_i^2$:
\begin{theorem}[Noiseless Point Sets with Outliers]\label{theorem:noiseless+outlier+c}
	Suppose there is no noise. Consider \eqref{eq:TLS-Q} with outliers $\bQ_{k^*+1},\dots,\bQ_\ell$ ($k^*<\ell$). Recall $\bw_0^*$ denotes the unit quaternion that represents the ground-truth rotation $\bR_0^*$. Let $\pw^*:=[\bw_0^*;\dots,\bw_0^*; 0;\dots; 0]\in\bbR^{4(\ell + 1)}$, where $\bw_0^*$ appeared $k^*+1$ times, and let $\pW^*:=\pw^*\trsp{(\pw^*)}$. Then we have:
	\begin{itemize}
		\item If $0<c_j^2<\lambda_{\textnormal{min}}(\bQ_{j})$,  $\forall j=k^*+1,\dots,\ell$, then \eqref{eq:SDP} is always tight, admitting $\pW^*$ as a global minimizer. 
		\item Suppose $c_j^2>\trsp{(\bw_0^*)}\bQ_j \bw_0^*$ for some $j\in\{ k^*+1,\dots,\ell\}$. Then $\pW^*$ is not a global minimizer of \eqref{eq:SDP}.
	\end{itemize}
\end{theorem}
\begin{proof}[Proof of Theorem \ref{theorem:noiseless+outlier+c}]
	For the first part, note that $c_j^2<\lambda_{\textnormal{min}}(\bQ_{j})$,  $\forall j=k^*+1,\dots,\ell$, so \eqref{eq:TLS-Q} rejects all outliers and preserves all inliers, $\bw_0^*$ globally minimizes \eqref{eq:TLS-Q}, and $\pw^*$ globally minimizes \eqref{eq:P-QCQP} with the minimum value $\sum_{j=k^*+1}^{\ell} c_j^2$. Let $[\hat{\pD}]_{0i}:= \frac{1}{2}(\bQ_i +c_i^2\bI_4)$ ($\forall i=1,\dots,k^*$),  $[\hat{\pD}]_{0j}:= \frac{1}{2}(\bQ_j - c_j^2\bI_4)$ ($\forall j=k^*+1,\dots,\ell$), and  $\hat{\mu} :=-\sum_{i=1}^{k^*}c_i^2$. With $\bQ_i\bw^*_0=0$, $\forall i=1,\dots,k^*$,  (\cref{lemma:rotation-quaternion}), one easily verifies conditions \eqref{eq:stationary2} and \eqref{eq:value-mu2} of \cref{prop:tightnesscondition} hold. Condition \eqref{eq:dualfeasibility2} is the same as
		\begin{equation*}
			\sum_{i=1}^{k^*} \Big(c_i^2\norm{\bz_0 - \bz_i}^2 + \trsp{\bz_{i}}  \bQ_i \bz_i \Big) + \sum_{j=k^*+1}^{\ell}\trsp{\bz_{j}}  (\bQ_j - c_j^2\bI_4) \bz_j \geq 0,\ \  \forall \bz_i\in\bbR^4,
		\end{equation*}
	which holds true because $\bQ_i$' are positive semidefinite as per \cref{lemma:rotation-quaternion} and $\bQ_{j}\succeq c_j^2\bI_4$. 
	
	For the second part, assume  $c_\ell^2>\trsp{(\bw_0^*)}\bQ_{\ell}\bw_0^*$, and we will show $\pw^*\trsp{(\pw^*)}$ does not minimize \eqref{eq:SDP}. For this we need to show that conditions \eqref{eq:stationary2}, \eqref{eq:dualfeasibility2}, and \eqref{eq:value-mu2} of \cref{prop:tightnesscondition} can not be simultaneously satisfied by $\bw_0^*$ and any $\hat{\mu}$ and any $\hat{\pD}$, where $\hat{\pD}$ is of the form \eqref{eq:D-constraint}. Assume  \eqref{eq:stationary2} and \eqref{eq:value-mu2} are true and we will prove \eqref{eq:dualfeasibility2} is false. In particular, \eqref{eq:value-mu2} implies $\hat{\mu}=-\sum_{i=1}^{k^*}c_i^2$. And  \eqref{eq:stationary2} implies $\big(2[\hat{\pD}]_{0\ell} + c_\ell^2\bI_4 -\ \bQ_\ell \big) \bw_0^* = 0$, which yields the inequality 
		\begin{equation*}
			2  \trsp{(\bw_0^*)}[\hat{\pD}]_{0\ell}\bw_0^* = \trsp{(\bw_0^*)}\bQ_{\ell}\bw_0^* - c_\ell^2 < 0.
		\end{equation*}
	Hence $2[\hat{\pD}]_{0\ell}$ must have a negative eigenvalue. Let $\pz=[\bz_0;\dots;\bz_\ell]\in\bbR^{4(\ell+1)}$ be such that $\bz_i=0$ for any $i\neq \ell$ and $\bz_{\ell}\in\bbS^3$ is an eigenvector of $2[\hat{\pD}]_{0\ell}$ corresponding to that negative eigenvalue. Then $\trsp{\pz}(\pQ-\hat{\mu} \pB-\hat{\pD})\pz=2\trsp{\bz_\ell}[\hat{\pD}]_{0\ell}\bz_\ell<0$, and condition \eqref{eq:dualfeasibility2} is not true. This finishes the proof of the second statement.
\end{proof}
\begin{remark}[Noiseless Point Sets with Random Outliers]\label{remark:random-outliers}
	If outlier $(\by_j,\bx_j)$ is randomly drawn from $\bbR^3\times \bbR^3$ according to some continuous probability distribution, then with probability $1$ we have $\|\by_j \|_2\neq \|\bx_j\|_2$ (Lemma 2 of \cite{Unnikrishnan-TIT18}), which implies $\lambda_{\textnormal{min}}(\bQ_{j})>0$. Thus, \eqref{eq:SDP} is always tight to such random outliers, if $c_j^2\to 0$. Note that this discussion is theoretical and does not apply to the case where $\big|\|\by_j \|_2- \|\bx_j\|_2 \big|$ is nonzero but is below machine accuracy, as $c_j^2$ can not be set even smaller (we only consider $c_j^2>0$).
\end{remark}
If the condition $c_j^2<\lambda_{\textnormal{min}}(\bQ_{j})$ of the first statement in \cref{theorem:noiseless+outlier+c} holds then $\bQ_{j}$ will always be rejected by \eqref{eq:TLS-Q} as an outlier. In fact, since $\lambda_{\textnormal{min}}(\bQ_{j})$ can be easily computed (\cref{lemma:rotation-quaternion}), in practice one usually throws away the point pairs $(\by_j,\bx_j)$'s for which $c_j^2<\lambda_{\textnormal{min}}(\bQ_{j})$ as a means of preprocessing (cf. \cite{Bustos-ICCV2015,Peng-CVPR2022}), and these point pairs might not enter into the semidefinite optimization. Our result here thus shows that \eqref{eq:SDP} can distinguish this type of ``simple'' outliers, as long as $c_j^2$ is chosen sufficiently small.

In the second statement of \cref{theorem:noiseless+outlier+c}, if $c_\ell^2>\trsp{(\bw_0^*)}\bQ_\ell \bw_0^*$ holds true for outlier $\bQ_\ell$, then \eqref{eq:TLS-Q} would attempt to minimize $\trsp{\bw_0} \bQ_\ell \bw_0  + \sum_{i=1}^{k^*} \trsp{\bw_0} \bQ_i \bw_0 $ over $\bw_0\in\bbS^3$ at least---an outlier showed up in the eigenvalue optimization---thus the global minimizer of \eqref{eq:TLS-Q} is unlikely to be $\bw_0^*$. This is why we do not expect $\pW^*$ to globally minimize \eqref{eq:SDP}.

Admittedly, \cref{theorem:noiseless+outlier+c} leaves a gap: What can we say about the tightness of \eqref{eq:SDP} if
\begin{equation}\label{eq:c-sandwich}
	\lambda_{\textnormal{min}}(\bQ_j) < c_j^2 < \trsp{(\bw_0^*)}\bQ_{j}\bw_0^* \text{\ \ ?}
\end{equation}
While our empirical observation suggests that $\pW^*$ does not globally minimize \eqref{eq:SDP} if \eqref{eq:c-sandwich} holds (with $k^*<\ell$), the analysis of this case  without further assumptions on outliers appears hard. The difficulty is that the outliers $\bQ_{j}$'s could be so adversarial that $\trsp{(\bw_0^*)}\bQ_{j}\bw_0^*$ is arbitrarily close\footnote{Alternatively, if $\trsp{(\bw_0^*)}\bQ_{j}\bw_0^*$ is small, then $\bQ_{j}$ might be treated as noisy data rather than an outlier. We consider such noisy case in \cref{subsection:noisy+outlier-free} (without outliers) and \cref{subsection:noisy+outliers} (with outliers).} to $0$ while $\lambda_{\textnormal{min}}(\bQ_j)=0$ for every $j>k^*$. Thus, the value of \cref{theorem:noiseless+outlier+c} is in that it shows that \eqref{eq:SDP} can only handle ``simple'' outliers that can be filtered out, and thus reveals a fundamental limit on the performance of \eqref{eq:SDP}.

Next, we consider the situation where outliers $\bQ_j$'s can not be simply removed by preprocessing, e.g., $\lambda_{\textnormal{min}}(\bQ_j)=0$. In particular, we assume the outliers are \textit{clustered} and show that the \eqref{eq:SDP} under investigation is even more vulnerable:
\begin{theorem}[Noiseless Point Sets with Clustered Outliers]\label{theorem:noiseless+outlier+c2} With the notation of \cref{theorem:noiseless+outlier+c}, further suppose outliers  $\bQ_{k^*+1},\dots,\bQ_{\ell}$ are ``clustered'' in the sense that
	\begin{equation}\label{eq:adversarial-outliers}
		\bQ_{k^*+1}\bw_0^{\textnormal{cl}}=\cdots =\bQ_{\ell}\bw_0^{\textnormal{cl}}=0
	\end{equation} 
	with $\bw_0^{\textnormal{cl}}\in\bbS^3$ some unit quaternion that is different from $\pm\bw_0^*$. If
	\begin{equation}\label{eq:condition-adversarial-outliers}
		1 - \frac{\sum_{j=k^*+1}^{\ell} c_j^2}{2 \sum_{i=1}^{k^*}c_i^2} <  \big| \trsp{(\bw_0^{\textnormal{cl}})} \bw_0^*  \big|,
	\end{equation}
	then $\pW^*$ does not globally minimize \eqref{eq:SDP}.
\end{theorem}
\begin{proof}[Proof of \cref{theorem:noiseless+outlier+c2}]
	Here we follow the proof idea in \cref{theorem:noiseless+outlier+c}. Let $\hat{\mu}:=-\sum_{i=1}^{k^*}c_i^2$, $\hat{\pD}$  of the form \eqref{eq:D-constraint}, and they satisfy \eqref{eq:stationary2} (and \eqref{eq:value-mu2}) with $\bw_0^*$. We will show that, given the clustered outliers which satisfy \eqref{eq:adversarial-outliers} and \eqref{eq:condition-adversarial-outliers}, the dual feasibility condition \eqref{eq:dualfeasibility2} is not true. Suppose without loss of generality $\trsp{(\bw_0^{\textnormal{cl}})} \bw_0^*>0$. With $\bz_1=\cdots=\bz_{k^*}=\bw_0^*$, and $\bz_{k^*+1}=\cdots=\bz_\ell = \bz_0=\bw_0^{\textnormal{cl}}$, \eqref{eq:dualfeasibility2} becomes
	\begin{equation*}
		\begin{split}
			&\ \ \sum_{i=1}^{k^*}c_i^2 \Big( \norm{ \bw_0^{\textnormal{cl}} }^2 + \norm{ \bw_0^* }^2 -  2 \trsp{(\bw_0^{\textnormal{cl}})} \bw_0^* \Big) +\sum_{j=k^*+1}^{\ell} \trsp{(\bw_0^{\textnormal{cl}})} \Big( \bQ_j - c_j^2\bI_4 \Big) \bw_0^{\textnormal{cl}} \geq 0 \\
			\Leftrightarrow & \ \  1 - \frac{\sum_{j=k^*+1}^{\ell} c_j^2}{2 \sum_{i=1}^{k^*}c_i^2} \geq    \trsp{(\bw_0^{\textnormal{cl}})} \bw_0^*,
		\end{split}
	\end{equation*}
	which violates assumption \eqref{eq:condition-adversarial-outliers}. This finishes the proof.
\end{proof}
The \textit{clustered} outliers of \cref{theorem:noiseless+outlier+c2} defined in the sense of \eqref{eq:adversarial-outliers} mean that the outlier pairs $(\by_j,\bx_j)$ ($j>k^*$) are related by the same 3D rotation $\bR_0^{\textnormal{cl}}$ that correspond to $\bw_0^{\textnormal{cl}}$, that is $\by_j=\bR_0^{\textnormal{cl}}\bx_j$, $\forall j>k^*$ (\cref{lemma:rotation-quaternion}). Clustered outliers can be thought of as a special type of \textit{adversarial} outliers, the latter usually used to study the robustness of algorithms in the worse case; it should be distinguished from data clustering \cite{Jain-CSUR1999,Elhamifar-PAMI2013}.

To understand condition \eqref{eq:condition-adversarial-outliers} of \cref{theorem:noiseless+outlier+c2}, consider a situation where all truncation parameters are equal, $c_1^2=\cdots=c_\ell^2$. Then \eqref{eq:condition-adversarial-outliers} simplifies to $1 - (\ell - k^*)/(2k^*) < \big| \trsp{(\bw_0^{\textnormal{cl}})} \bw_0^*  \big|$; also note that $\big| \trsp{(\bw_0^{\textnormal{cl}})} \bw_0^*  \big|\in[0,1)$. Thus, if $\ell-k^*>2k^*$, then \eqref{eq:condition-adversarial-outliers} always holds, and so $\pW^*$ never globally minimizes \eqref{eq:SDP}, which is forgivable as in this case $\bw_0^*$ neither globally minimizes \eqref{eq:TLS-Q}. However, even if the number of outliers is only half the number of inliers, i.e., $\ell-k^*=k^*/2$, \cref{theorem:noiseless+outlier+c2} implies that $\pW^*$ would still fail to globally minimize \eqref{eq:SDP} as long as $\big| \trsp{(\bw_0^{\textnormal{cl}})} \bw_0^*  \big|> 3/4$, but $\bw_0^*$ would  \textit{in general} globally minimize \eqref{eq:TLS-Q} with suitable $c_j^2$ (cf. \cite{Yang-T-R2021}). In this sense, \eqref{eq:SDP} is \textit{strictly} less robust to outliers than \eqref{eq:TLS-Q}.

Finally, we note that \cref{theorem:noiseless+outlier+c2} might be overly pessimistic. In fact, experiments show that \eqref{eq:SDP} is robust to $40\%$-$50\%$ outliers $(\by_j,\bx_j)$'s, where $\by_j$ and $\bx_j$ are sampled uniformly at random from $\bbS^2$ (so $\lambda_{\textnormal{min}}(\bQ_j)=0$ by \cref{lemma:rotation-quaternion}). Two factors account for this empirically better behavior: i) Such random outliers are less adversarial than clustered ones, ii) the extra projection step that converts the global minimizer of \eqref{eq:SDP} to a unit quaternion alleviates  to some extent the issue of $\pW^*$ not minimizing \eqref{eq:SDP}. In retrospect, there are two  downsides in our analysis of Theorems \ref{theorem:noiseless+outlier+c} and \ref{theorem:noiseless+outlier+c2}: (1) We have not taken such extra projection step into account, and (2) we only showed that $\pW^*$ might not minimize \eqref{eq:SDP} but have not proved how far the global minimizers of \eqref{eq:SDP} can be from $\pW^*$. 

\subsection{The Noisy + Outlier-Free Case}\label{subsection:noisy+outlier-free}
The noisy case, even without outliers, is more difficult to penetrate than previous cases. A general reason for this is that the global minimizers of \eqref{eq:P-QCQP} and \eqref{eq:SDP} are now complicated functions of noise. One might wonder whether \cref{theorem:noiseless+outlier-free} can be extended to the noisy + outlier-free case using some continuity argument. In fact, \cite{Cifuentes-MP2021} shows that, under certain conditions, if the noiseless version of the Schor relaxation of some QCQP is tight, then its noisy version is also tight. While this result is quite general, its conditions are abstract and, at least for non-experts, hard to verify. To our case, \cite{Cifuentes-MP2021} is not directly applicable, as in our problem the truncation parameters $c_i^2$ also have impacts on the tightness, and the approach of \cite{Cifuentes-MP2021} does not model, and thus could not control the values of $c_i^2$. Instead, our analysis must take both $c_i^2$ and noise into account.

We begin our analysis by decomposing the corrupted data matrix $\bQ_i$ of \eqref{eq:TLS-Q} into 
the pure data part $\bP_i$ and noise part $\bE_i+\|\beps_i \|_2^2\bI_4$:
\begin{lemma}\label{lemma:rotation-quaternion2}
	Let $\bR_0\in\SO(3)$. If $(\by_i,\bx_i)$ is an inlier that obeys \eqref{eq:RRS-model}, then  we have
	\begin{equation}\label{eq:Qtilde=data+noise}
		\norm{\by_i-\bR_0\bx_i}^2=\trsp{\bw_0} \bQ_i \bw_0,\ \ \ \bQ_i =\bP_i +\bE_i+\norm{\beps_i}^2\bI_4,
	\end{equation}
	where $\bw_0\in\bbS^3$ is the unit quaternion representation of $\bR_0$, and $\bP_i$ and $\bE_i$ are $4\times 4$ symmetric matrices that repsectively satisfy the following properties:
	\begin{itemize}
		\item $\bP_i$ is positive semidefinite with its entries depending on $\by_i$ and $\bx_i$, and it has two different eigenvalues $4\|\bx_i\|_2^2$ and $0$, each of multiplicity $2$. The ground-truth unit quaternion $\bw_0^*$ is an eigenvector of $\bP_i$ corresponding to eigenvalue $0$, i.e., $\bP_i\bw_0^*=0$. In particular, we have $\bQ_i\bw_0^*=\bP_i\bw_0^*=0$ in the noiseless case.
		\item $\bE_i$ has entries depending on $\by_i$, $\bx_i$, and noise  $\beps_i$, and  has two different eigenvalues $2\trsp{\beps_i}\bR_0^*\bx_i+2\|\beps_i \|_2 \| \bx_i \|_2$ and $2\trsp{\beps_i}\bR_0^*\bx_i- 2\|\beps_i \|_2 \| \bx_i \|_2$, each of multiplicity $2$. We have $\trsp{\bw_0}\bE_i\bw_0 = 2\trsp{\beps_i}(\bR_0^*\bx_i - \bR_0\bx_i)$ and in particular $\trsp{(\bw_0^*)}\bE_i\bw_0^*=0$. 
	\end{itemize}
\end{lemma}
\begin{proof}[Proof of \cref{lemma:rotation-quaternion2}]
	We start with the following equality:
	\begin{equation*}
		\norm{\by_i-\bR_0\bx_i}^2=\norm{\beps_i+\bR_0^*\bx_i-\bR_0\bx_i}^2
		=\norm{\bR_0^*\bx_i-\bR\bx_i}^2 + \norm{\beps_i}^2 + 2\trsp{\beps_i}(\bR_0^*\bx_i - \bR_0\bx_i) 
	\end{equation*}
	From the proof of \cref{lemma:rotation-quaternion} and in particular \eqref{eq:Q(y,x)}, we see that the first term $\|\bR_0^*\bx_i-\bR_0\bx_i\|_2^2$ is equal to $\trsp{\bw_0}\bP_i\bw_0$, where $\bP_i:=Q(\bR_0^*\bx_i,\bx_i)$ is a $4\times 4$ positive semidefinite matrix with eigenvalues $4\|\bx_i\|_2^2, 4\|\bx_i\|^2, 0,0$, and $\bw_0$ is the unit quaternion that represents $\bR_0$. In particular, we have $\trsp{(\bw_0^*)}\bP_i\bw_0^*=\|\bR_0^*\bx_i-\bR_0^*\bx_i\|_2^2=0$, which implies $\bP_i\bw_0^*=0$. 
	
	On the other hand, from the proof of \cref{lemma:rotation-quaternion} and in particular \eqref{eq:def-Qi}, we see the equality $-2\trsp{\beps_i}\bR_0\bx_i=-2\trsp{\bw_0} \ U(\beps_i,\bx_i)\ \bw_0$ where $-2U(\beps_i,\bx_i)$ is symmetric and has two different eigenvalues $2\|\beps_i\|_2\|\bx_i\|_2$ and $-2\|\beps_i\|_2\|\bx_i\|_2$, each of multiplicity $2$. Then $\bE_i:=2(\trsp{\beps_i}\bR_0^*\bx_i)\bI_4 -2 U(\beps_i,\bx_i)$ satisfies the claimed properties.
\end{proof}

Our investigation into the noisy case is tightly related to 
the \textit{eigengap} $\zeta$, defined as the ratio between the second smallest eigenvalue $\lambda_{\textnormal{min2}}(\cdot )$ of $\sum_{i=1}^{\ell} 	\bQ_i$ and its minimum eigenvalue:
\begin{equation}\label{eq:eigengap}
	\zeta:=\frac{\lambda_{\textnormal{min2}}\big(  \sum_{i=1}^{\ell} 	\bQ_i \big)}{\lambda_{\textnormal{min}}\big(  \sum_{i=1}^{\ell} 	\bQ_i \big)}
\end{equation}
In analysis of eigenvalue algorithms (cf. \cite{Chen-FTML2021}), the eigengap is typically defined as the \textit{difference} between two consecutive eigenvalues of some matrix. Our eigengap \eqref{eq:eigengap} that takes the \textit{division} of two smallest eigenvalues is not standard, but it will be convenient for our purpose. Also note that, $\bQ_i$ is a perturbed version of $\bP_i$ by noise $\beps_i$ (\cref{lemma:rotation-quaternion2}), so $\lambda_{\textnormal{min}}\big(  \sum_{i=1}^{\ell} 	\bQ_i \big)$ is in general nonzero, while it indeed approaches zero if $\beps_i\to 0$. Clearly $\zeta\geq 1$. Moreover, we have the following immediate observation:
\begin{remark}\label{remark:eigengap}
	If \eqref{eq:TLS-Q} has a unique solution and if $c_i^2> \trsp{\hat{\bw}_0} \bQ_i \hat{\bw}_0$ ($\forall i$), then $\zeta>1$.	
\end{remark}

We are now ready to state the following result:
\begin{theorem}[Noisy and Outlier-Free Point Sets]\label{theorem:noisy+outlier-free2}
	Consider \eqref{eq:TLS-Q} with noisy inliers $\{\bQ_i\}_{i=1}^\ell$ and $\hat{\bw}_0\in\bbS^3$ its global minimizer. Let $\hat{\pw}:=[\hat{\bw}_0;\dots;\hat{\bw}_0]\in\bbR^{4(\ell+1)}$. Suppose $\zeta \geq \ell/(\ell-1)$. \eqref{eq:SDP} is tight as long as 
	\begin{equation}
		\begin{split}
			&c_i^2 > \trsp{\hat{\bw}_0} \bQ_i \hat{\bw}_0 + \norm{\bQ_i \hat{\bw}_0} + \frac{|d_i| + d_i}{2}, \ \ \ \forall i=1,\dots,\ell    \\
			\textnormal{with \ } & d_i := \frac{\sum_{i=1}^{\ell} \trsp{\hat{\bw}_0}\bQ_i\hat{\bw}_0}{\ell} - \trsp{\hat{\bw}_0} \bQ_i \hat{\bw}_0 + \frac{\lambda_{\textnormal{max}}\Big(\sum_{j\neq i} \big(  \bQ_{i}  -  \bQ_{j} \big)\Big)}{\zeta(\ell - 1)}.
		\end{split} \label{eq:noisy+outlier-free2-condition}
	\end{equation}
	Moreover, the angle $\hat{\tau}^*_0$ between $\hat{\bw}_0$ and the ground-truth unit quaternion $\bw^*_0\in \bbS^3$ grows proportionally with the magnitude of noise $\beps_i$:
	\begin{equation}\label{eq:dist-eigvector2gt}
		\sin^2(\hat{\tau}^*_0) \leq \frac{4\sum_{i=1}^{\ell}\norm{\beps_i} \norm{\bx_i}}{\lambda_{\textnormal{min2}}\big(\sum_{i=1}^{\ell} \bP_i\big)},  \ \ \ \ \  \sin^2(\hat{\tau}^*_0):= 1 - (\trsp{\hat{\bw}_0}\bw_0^*)^2
	\end{equation}
	In \eqref{eq:dist-eigvector2gt}, each $\bP_i$ corresponds to the ``pure data'' part of $\bQ_i$ that satisfies $\bP_i\bw_0^*=0,\bP_i\succeq0$ (\cref{lemma:rotation-quaternion2}), and $\lambda_{\textnormal{min2}}\big(\cdot)$ denotes the second smallest eigenvalue of a matrix.
\end{theorem}
\begin{proof}[Proof of Theorem \ref{theorem:noisy+outlier-free2}]
	Since \eqref{eq:noisy+outlier-free2-condition} implies $c_i^2 > \trsp{\hat{\bw}_0} \bQ_i \hat{\bw}_0$, \eqref{eq:TLS-Q} preserves all inliers with minimum $\sum_{i=1}^{\ell} \trsp{\hat{\bw}_0}\bQ_i\hat{\bw}_0$. So $\hat{\mu}:=\sum_{i=1}^{\ell}\big( \trsp{\hat{\bw}_0}\bQ_i\hat{\bw}_0 - c_i^2 \big)$ satisfies \eqref{eq:value-mu2} of \cref{prop:tightnesscondition}. 
	
	Let $\widehat{\bV}:=[\widehat{\bV}_0,\ \ \hat{\bw}_0]\in\bbR^{4\times 4}$ form an orthonormal basis of $\bbR^4$; $\widehat{\bV}_0\in\bbR^{4\times 3}$ satisfies $\trsp{\widehat{\bV}_0} \hat{\bw}_0 = 0$ and  $\trsp{\widehat{\bV}_0} \widehat{\bV}_0 = \bI_3$. Let  $\hat{\pD}$ be of the form \eqref{eq:D-constraint} with each $[\hat{\pD}]_{0i}$ satisfying 
	\begin{equation}\label{eq:D0i-Si}
		[\hat{\pD}]_{0i} := \widehat{\bS}_{i} - \frac{1}{2}\Big(\bQ_i -c_i^2\bI_4\Big),\ \ \ \widehat{\bS}_{i} := \widehat{\bV} \begin{bmatrix}
			\widehat{\bT}_i  & 0 \\
			0 & 0
		\end{bmatrix} \trsp{\widehat{\bV}}
	\end{equation}
	where $\widehat{\bT}_i$ is a  $3\times 3$ symmetric matrix defined as
	\begin{equation}\label{eq:Ti}
		\widehat{\bT}_i :=\frac{\zeta-\frac{\ell}{\ell-1}}{\zeta}\trsp{\widehat{\bV}_0} \bQ_{i} \widehat{\bV}_0 + \frac{\sum_{j= 1}^\ell \trsp{\widehat{\bV}_0} \bQ_{j} \widehat{\bV}_0}{\zeta(\ell - 1)} - \frac{\Big( \sum_{i=1}^{\ell} \trsp{\hat{\bw}_0}\bQ_i\hat{\bw}_0 \Big) \bI_3}{\ell}.
	\end{equation}
	Clearly  $\widehat{\bS}_{i}\hat{\bw}_0=0$, with which one verifies that \eqref{eq:stationary2} of \cref{prop:tightnesscondition} is fulfilled.
	
	It remains to prove \eqref{eq:dualfeasibility2} holds true. For this we need the following lemma.  \cref{lemma:property-Si} is a bit technical, so we put its proof into the appendix.
	\begin{lemma}\label{lemma:property-Si}
		Let $[\hat{\pD}]_{0i}, \widehat{\bS}_{i}$ and $\widehat{\bT}_{i}$ be defined in \eqref{eq:D0i-Si} and \eqref{eq:Ti} respectively. Suppose that the assumptions of \cref{theorem:noisy+outlier-free2}, namely $\zeta \geq \ell/(\ell-1)$ and \eqref{eq:noisy+outlier-free2-condition}, hold. Then
		\begin{equation*}
			\begin{split}
				\sum_{i=1}^{\ell} \widehat{\bS}_i &= \sum_{i=1}^{\ell} \big( \bQ_i  -  \trsp{\hat{\bw}_0}\bQ_i\hat{\bw}_0 \bI_4  \big) \\
				\widehat{\bS}_{i}&\succeq 0, \ \ \ \ \ \ \forall i=1,\dots,\ell  \\
				\widehat{\bS}_{i} + c_i^2\bI_4 -  \bQ_i&\succ 0, \ \ \ \ \ \ \forall i=1,\dots,\ell 
			\end{split}
		\end{equation*}
	\end{lemma}
	With \cref{lemma:property-Si} and definitions of $[\hat{\pD}]_{0i}$ \eqref{eq:D0i-Si}, we can now write \eqref{eq:dualfeasibility2} as ($\forall \bz_i\in\bbR^4$)
	\begin{equation*}
		\begin{split}
			&\sum_{i=1}^{\ell}  \Big( \trsp{\bz_0} \big( c_i^2 - \trsp{\hat{\bw}_0}\bQ_i\hat{\bw}_0   \big) \bz_0+ 
			\trsp{\bz_i} \big(2\widehat{\bS}_{i} + c_i^2\bI_4 -  \bQ_i\big)\bz_i - 
			2\trsp{\bz_0} \big( \widehat{\bS}_{i} + c_i^2\bI_4 -  \bQ_i  \big)\bz_i  \Big) \geq 0 \\
			\Leftrightarrow & \sum_{i=1}^{\ell}  \Big( \trsp{\bz_0} \big( \widehat{\bS}_{i} + c_i^2\bI_4 -  \bQ_i   \big) \bz_0+ 
			\trsp{\bz_i} \big(2\widehat{\bS}_{i} + c_i^2\bI_4 -  \bQ_i\big)\bz_i - 
			2\trsp{\bz_0} \big( \widehat{\bS}_{i} + c_i^2\bI_4 -  \bQ_i  \big)\bz_i  \Big) \geq 0 \\
			\Leftrightarrow& \sum_{i=1}^{\ell}	 \trsp{(\bz_i - \bz_0)}\Big(\widehat{\bS}_{i} + c_i^2\bI_4  - \bQ_i \Big)(\bz_i - \bz_0) + \sum_{i=1}^{\ell}   \trsp{\bz_i} \widehat{\bS}_{i} \bz_i \geq 0
		\end{split}
	\end{equation*}
	which holds, as $\widehat{\bS}_{i} + c_i^2\bI_4  - \bQ_i$ and $\widehat{\bS}_{i}$ are both positive semidefinite for every $i=1,\dots,\ell$ (\cref{lemma:property-Si}). We have thus verified all conditions of \cref{prop:tightnesscondition}, thus \eqref{eq:SDP} is tight.
	
	It now remains to prove the error bound \eqref{eq:dist-eigvector2gt}. With \cref{lemma:rotation-quaternion2} we write $\bQ_i =\bP_i +\bE_i+\|\beps_i\|_2^2\bI_4$, and we know that $\hat{\bw}_0$ is an eigenvector of $\sum_{i=1}^{\ell} (\bP_i +\bE_i)$ corresponding to its minimum eigenvalue. Hence, with the optimality of $\hat{\bw}_0$ and \cref{lemma:rotation-quaternion2}, we have
	\begin{equation*}
		\begin{split}
			& \sum_{i=1}^{\ell} \trsp{\hat{\bw}_0} (\bP_i +\bE_i ) \hat{\bw}_0 \leq \sum_{i=1}^{\ell} \trsp{(\bw_0^*)} (\bP_i +\bE_i ) \bw_0^* = 0 \\
			\Rightarrow & \sum_{i=1}^{\ell} \trsp{\hat{\bw}_0} \bP_i  \hat{\bw}_0 \leq - \sum_{i=1}^{\ell} \trsp{\hat{\bw}_0}\bE_i  \hat{\bw}_0 = - \sum_{i=1}^{\ell} 2\trsp{\beps_i}(\bR_0^* - \widehat{\bR}_0 ) \bx_i \leq  4\sum_{i=1}^{\ell}\norm{\beps_i} \norm{\bx_i} \nonumber
		\end{split}
	\end{equation*}
	Here we recall that $\widehat{\bR}_0$ is the global minimizer of \eqref{eq:TLS-R} that corresponds to $\hat{\bw}_0$, and used \cref{lemma:rotation-quaternion2}. Since $\sum_{i=1}^{\ell}\bP_i$ is positive semidefinite with  $\bP_i \bw_0^*=0$, from the eigen decomposition of $\sum_{i=1}^{\ell}\bP_i$ we can easily obtain the inequality
	\begin{equation*}
		\big( 1- (\trsp{\hat{\bw}_0}\bw_0^*)^2 \big) \lambda_{\textnormal{min2}}\Big(\sum_{i=1}^{\ell} \bP_i\Big) \leq \sum_{i=1}^{\ell} \trsp{\hat{\bw}_0} \bP_i  \hat{\bw}_0,
	\end{equation*}
	where $\lambda_{\textnormal{min2}}(\cdot)$ denotes the second smallest eigenvalue of a matrix. The proof is complete. 
\end{proof}

\begin{figure}[t]
	\centering
	\subfloat[]{
\begin{tikzpicture}[x=1pt,y=1pt]
\definecolor{fillColor}{RGB}{255,255,255}
\path[use as bounding box,fill=fillColor,fill opacity=0.00] (0,0) rectangle (122.86, 93.95);
\begin{scope}
\path[clip] (  0.00,  0.00) rectangle (122.86, 93.95);
\definecolor{drawColor}{RGB}{255,255,255}
\definecolor{fillColor}{RGB}{255,255,255}

\path[draw=drawColor,line width= 0.6pt,line join=round,line cap=round,fill=fillColor] ( -0.00,  0.00) rectangle (122.86, 93.95);
\end{scope}
\begin{scope}
\path[clip] ( 29.39, 27.67) rectangle ( 99.40, 88.45);
\definecolor{fillColor}{gray}{0.92}

\path[fill=fillColor] ( 29.39, 27.67) rectangle ( 99.40, 88.45);
\definecolor{drawColor}{RGB}{160,32,240}

\path[draw=drawColor,line width= 0.6pt,line join=round] ( 32.58, 75.36) --
	( 39.65, 76.05) --
	( 46.72, 76.18) --
	( 53.79, 74.25) --
	( 60.86, 74.96) --
	( 67.93, 77.31) --
	( 75.01, 75.95) --
	( 82.08, 77.41) --
	( 89.15, 76.58) --
	( 96.22, 75.03);

\path[draw=drawColor,line width= 0.4pt,line join=round,line cap=round] ( 31.15, 73.93) rectangle ( 34.00, 76.79);

\path[draw=drawColor,line width= 0.4pt,line join=round,line cap=round] ( 38.22, 74.63) rectangle ( 41.07, 77.48);

\path[draw=drawColor,line width= 0.4pt,line join=round,line cap=round] ( 45.29, 74.75) rectangle ( 48.15, 77.61);

\path[draw=drawColor,line width= 0.4pt,line join=round,line cap=round] ( 52.36, 72.82) rectangle ( 55.22, 75.67);

\path[draw=drawColor,line width= 0.4pt,line join=round,line cap=round] ( 59.44, 73.53) rectangle ( 62.29, 76.39);

\path[draw=drawColor,line width= 0.4pt,line join=round,line cap=round] ( 66.51, 75.88) rectangle ( 69.36, 78.74);

\path[draw=drawColor,line width= 0.4pt,line join=round,line cap=round] ( 73.58, 74.52) rectangle ( 76.43, 77.38);

\path[draw=drawColor,line width= 0.4pt,line join=round,line cap=round] ( 80.65, 75.98) rectangle ( 83.50, 78.83);

\path[draw=drawColor,line width= 0.4pt,line join=round,line cap=round] ( 87.72, 75.15) rectangle ( 90.58, 78.00);

\path[draw=drawColor,line width= 0.4pt,line join=round,line cap=round] ( 94.80, 73.61) rectangle ( 97.65, 76.46);
\definecolor{drawColor}{RGB}{0,0,0}

\path[draw=drawColor,line width= 0.6pt,line join=round] ( 32.58, 30.44) --
	( 39.65, 32.25) --
	( 46.72, 34.99) --
	( 53.79, 38.75) --
	( 60.86, 44.20) --
	( 67.93, 50.62) --
	( 75.01, 58.19) --
	( 82.08, 67.11) --
	( 89.15, 76.48) --
	( 96.22, 85.69);

\path[draw=drawColor,line width= 0.4pt,line join=round,line cap=round] ( 32.58, 32.66) --
	( 34.50, 29.33) --
	( 30.65, 29.33) --
	( 32.58, 32.66);

\path[draw=drawColor,line width= 0.4pt,line join=round,line cap=round] ( 39.65, 34.47) --
	( 41.57, 31.14) --
	( 37.73, 31.14) --
	( 39.65, 34.47);

\path[draw=drawColor,line width= 0.4pt,line join=round,line cap=round] ( 46.72, 37.21) --
	( 48.64, 33.88) --
	( 44.80, 33.88) --
	( 46.72, 37.21);

\path[draw=drawColor,line width= 0.4pt,line join=round,line cap=round] ( 53.79, 40.97) --
	( 55.71, 37.64) --
	( 51.87, 37.64) --
	( 53.79, 40.97);

\path[draw=drawColor,line width= 0.4pt,line join=round,line cap=round] ( 60.86, 46.42) --
	( 62.78, 43.09) --
	( 58.94, 43.09) --
	( 60.86, 46.42);

\path[draw=drawColor,line width= 0.4pt,line join=round,line cap=round] ( 67.93, 52.84) --
	( 69.86, 49.51) --
	( 66.01, 49.51) --
	( 67.93, 52.84);

\path[draw=drawColor,line width= 0.4pt,line join=round,line cap=round] ( 75.01, 60.40) --
	( 76.93, 57.08) --
	( 73.09, 57.08) --
	( 75.01, 60.40);

\path[draw=drawColor,line width= 0.4pt,line join=round,line cap=round] ( 82.08, 69.33) --
	( 84.00, 66.00) --
	( 80.16, 66.00) --
	( 82.08, 69.33);

\path[draw=drawColor,line width= 0.4pt,line join=round,line cap=round] ( 89.15, 78.70) --
	( 91.07, 75.38) --
	( 87.23, 75.38) --
	( 89.15, 78.70);

\path[draw=drawColor,line width= 0.4pt,line join=round,line cap=round] ( 96.22, 87.91) --
	( 98.14, 84.58) --
	( 94.30, 84.58) --
	( 96.22, 87.91);
\end{scope}
\begin{scope}
\path[clip] (  0.00,  0.00) rectangle (122.86, 93.95);
\definecolor{drawColor}{RGB}{160,32,240}

\node[text=drawColor,anchor=base east,inner sep=0pt, outer sep=0pt, scale=  0.59] at ( 24.44, 27.83) {0};

\node[text=drawColor,anchor=base east,inner sep=0pt, outer sep=0pt, scale=  0.59] at ( 24.44, 43.84) {250};

\node[text=drawColor,anchor=base east,inner sep=0pt, outer sep=0pt, scale=  0.59] at ( 24.44, 59.85) {500};

\node[text=drawColor,anchor=base east,inner sep=0pt, outer sep=0pt, scale=  0.59] at ( 24.44, 75.86) {750};
\end{scope}
\begin{scope}
\path[clip] (  0.00,  0.00) rectangle (122.86, 93.95);
\definecolor{drawColor}{gray}{0.20}

\path[draw=drawColor,line width= 0.6pt,line join=round] ( 26.64, 29.88) --
	( 29.39, 29.88);

\path[draw=drawColor,line width= 0.6pt,line join=round] ( 26.64, 45.89) --
	( 29.39, 45.89);

\path[draw=drawColor,line width= 0.6pt,line join=round] ( 26.64, 61.90) --
	( 29.39, 61.90);

\path[draw=drawColor,line width= 0.6pt,line join=round] ( 26.64, 77.91) --
	( 29.39, 77.91);
\end{scope}
\begin{scope}
\path[clip] (  0.00,  0.00) rectangle (122.86, 93.95);
\definecolor{drawColor}{gray}{0.20}

\path[draw=drawColor,line width= 0.6pt,line join=round] ( 99.40, 29.86) --
	(102.15, 29.86);

\path[draw=drawColor,line width= 0.6pt,line join=round] ( 99.40, 49.09) --
	(102.15, 49.09);

\path[draw=drawColor,line width= 0.6pt,line join=round] ( 99.40, 68.31) --
	(102.15, 68.31);

\path[draw=drawColor,line width= 0.6pt,line join=round] ( 99.40, 87.54) --
	(102.15, 87.54);
\end{scope}
\begin{scope}
\path[clip] (  0.00,  0.00) rectangle (122.86, 93.95);
\definecolor{drawColor}{RGB}{0,0,0}

\node[text=drawColor,anchor=base west,inner sep=0pt, outer sep=0pt, scale=  0.59] at (104.35, 27.82) {0};

\node[text=drawColor,anchor=base west,inner sep=0pt, outer sep=0pt, scale=  0.59] at (104.35, 47.04) {1};

\node[text=drawColor,anchor=base west,inner sep=0pt, outer sep=0pt, scale=  0.59] at (104.35, 66.27) {2};

\node[text=drawColor,anchor=base west,inner sep=0pt, outer sep=0pt, scale=  0.59] at (104.35, 85.49) {3};
\end{scope}
\begin{scope}
\path[clip] (  0.00,  0.00) rectangle (122.86, 93.95);
\definecolor{drawColor}{gray}{0.20}

\path[draw=drawColor,line width= 0.6pt,line join=round] ( 32.58, 24.92) --
	( 32.58, 27.67);

\path[draw=drawColor,line width= 0.6pt,line join=round] ( 53.79, 24.92) --
	( 53.79, 27.67);

\path[draw=drawColor,line width= 0.6pt,line join=round] ( 75.01, 24.92) --
	( 75.01, 27.67);

\path[draw=drawColor,line width= 0.6pt,line join=round] ( 96.22, 24.92) --
	( 96.22, 27.67);
\end{scope}
\begin{scope}
\path[clip] (  0.00,  0.00) rectangle (122.86, 93.95);
\definecolor{drawColor}{RGB}{0,0,0}

\node[text=drawColor,anchor=base,inner sep=0pt, outer sep=0pt, scale=  0.79] at ( 32.58, 17.27) {$1\%$};

\node[text=drawColor,anchor=base,inner sep=0pt, outer sep=0pt, scale=  0.79] at ( 53.79, 17.27) {$4\%$};

\node[text=drawColor,anchor=base,inner sep=0pt, outer sep=0pt, scale=  0.79] at ( 75.01, 17.27) {$7\%$};

\node[text=drawColor,anchor=base,inner sep=0pt, outer sep=0pt, scale=  0.79] at ( 96.22, 17.27) {$10\%$};
\end{scope}
\begin{scope}
\path[clip] (  0.00,  0.00) rectangle (122.86, 93.95);
\definecolor{drawColor}{RGB}{0,0,0}

\node[text=drawColor,anchor=base,inner sep=0pt, outer sep=0pt, scale=  0.85] at ( 64.40,  7.15) {Noise Level ($\sigma$)};
\end{scope}
\begin{scope}
\path[clip] (  0.00,  0.00) rectangle (122.86, 93.95);
\definecolor{drawColor}{RGB}{160,32,240}

\node[text=drawColor,rotate= 90.00,anchor=base,inner sep=0pt, outer sep=0pt, scale=  0.82] at ( 11.18, 58.06) {$\lambda_{\textnormal{min2}}\big(\sum_{i=1}^{\ell} Q_i\big)$};
\end{scope}
\begin{scope}
\path[clip] (  0.00,  0.00) rectangle (122.86, 93.95);
\definecolor{drawColor}{RGB}{0,0,0}

\node[text=drawColor,rotate=-90.00,anchor=base,inner sep=0pt, outer sep=0pt, scale=  0.82] at (110.07, 58.06) {$\lambda_{\textnormal{min}}\big(\sum_{i=1}^{\ell}Q_i\big)$};
\end{scope}
\end{tikzpicture} \label{fig:eigengap} }
	\subfloat[]{
\begin{tikzpicture}[x=1pt,y=1pt]
\definecolor{fillColor}{RGB}{255,255,255}
\path[use as bounding box,fill=fillColor,fill opacity=0.00] (0,0) rectangle (102.62, 93.95);
\begin{scope}
\path[clip] (  0.00,  0.00) rectangle (102.62, 93.95);
\definecolor{drawColor}{RGB}{255,255,255}
\definecolor{fillColor}{RGB}{255,255,255}

\path[draw=drawColor,line width= 0.6pt,line join=round,line cap=round,fill=fillColor] ( -0.00,  0.00) rectangle (102.62, 93.95);
\end{scope}
\begin{scope}
\path[clip] ( 25.07, 27.67) rectangle ( 97.12, 88.45);
\definecolor{fillColor}{gray}{0.92}

\path[fill=fillColor] ( 25.07, 27.67) rectangle ( 97.12, 88.45);
\definecolor{drawColor}{RGB}{0,0,0}

\path[draw=drawColor,line width= 0.6pt,line join=round] ( 28.34, 30.44) --
	( 35.62, 40.70) --
	( 42.90, 45.87) --
	( 50.18, 50.59) --
	( 57.46, 53.19) --
	( 64.73, 56.49) --
	( 72.01, 58.35) --
	( 79.29, 60.09) --
	( 86.57, 61.84) --
	( 93.85, 62.83);

\path[draw=drawColor,line width= 0.4pt,line join=round,line cap=round] ( 26.92, 29.01) rectangle ( 29.77, 31.86);

\path[draw=drawColor,line width= 0.4pt,line join=round,line cap=round] ( 34.19, 39.27) rectangle ( 37.05, 42.12);

\path[draw=drawColor,line width= 0.4pt,line join=round,line cap=round] ( 41.47, 44.44) rectangle ( 44.33, 47.29);

\path[draw=drawColor,line width= 0.4pt,line join=round,line cap=round] ( 48.75, 49.16) rectangle ( 51.60, 52.02);

\path[draw=drawColor,line width= 0.4pt,line join=round,line cap=round] ( 56.03, 51.76) rectangle ( 58.88, 54.61);

\path[draw=drawColor,line width= 0.4pt,line join=round,line cap=round] ( 63.31, 55.06) rectangle ( 66.16, 57.92);

\path[draw=drawColor,line width= 0.4pt,line join=round,line cap=round] ( 70.59, 56.92) rectangle ( 73.44, 59.78);

\path[draw=drawColor,line width= 0.4pt,line join=round,line cap=round] ( 77.86, 58.67) rectangle ( 80.72, 61.52);

\path[draw=drawColor,line width= 0.4pt,line join=round,line cap=round] ( 85.14, 60.41) rectangle ( 88.00, 63.26);

\path[draw=drawColor,line width= 0.4pt,line join=round,line cap=round] ( 92.42, 61.40) rectangle ( 95.27, 64.26);
\definecolor{drawColor}{RGB}{160,32,240}

\path[draw=drawColor,line width= 0.6pt,line join=round] ( 28.34, 68.87) --
	( 35.62, 73.78) --
	( 42.90, 76.33) --
	( 50.18, 78.64) --
	( 57.46, 79.92) --
	( 64.73, 82.01) --
	( 72.01, 82.93) --
	( 79.29, 84.06) --
	( 86.57, 84.57) --
	( 93.85, 85.69);

\path[draw=drawColor,line width= 0.4pt,line join=round,line cap=round] ( 26.92, 67.44) rectangle ( 29.77, 70.30);

\path[draw=drawColor,line width= 0.4pt,line join=round,line cap=round] ( 34.19, 72.35) rectangle ( 37.05, 75.21);

\path[draw=drawColor,line width= 0.4pt,line join=round,line cap=round] ( 41.47, 74.91) rectangle ( 44.33, 77.76);

\path[draw=drawColor,line width= 0.4pt,line join=round,line cap=round] ( 48.75, 77.22) rectangle ( 51.60, 80.07);

\path[draw=drawColor,line width= 0.4pt,line join=round,line cap=round] ( 56.03, 78.50) rectangle ( 58.88, 81.35);

\path[draw=drawColor,line width= 0.4pt,line join=round,line cap=round] ( 63.31, 80.58) rectangle ( 66.16, 83.43);

\path[draw=drawColor,line width= 0.4pt,line join=round,line cap=round] ( 70.59, 81.50) rectangle ( 73.44, 84.35);

\path[draw=drawColor,line width= 0.4pt,line join=round,line cap=round] ( 77.86, 82.63) rectangle ( 80.72, 85.48);

\path[draw=drawColor,line width= 0.4pt,line join=round,line cap=round] ( 85.14, 83.14) rectangle ( 88.00, 86.00);

\path[draw=drawColor,line width= 0.4pt,line join=round,line cap=round] ( 92.42, 84.26) rectangle ( 95.27, 87.12);
\end{scope}
\begin{scope}
\path[clip] (  0.00,  0.00) rectangle (102.62, 93.95);
\definecolor{drawColor}{RGB}{0,0,0}

\node[text=drawColor,anchor=base east,inner sep=0pt, outer sep=0pt, scale=  0.73] at ( 20.12, 35.39) {5E-4};

\node[text=drawColor,anchor=base east,inner sep=0pt, outer sep=0pt, scale=  0.73] at ( 20.12, 51.60) {5E-3};

\node[text=drawColor,anchor=base east,inner sep=0pt, outer sep=0pt, scale=  0.73] at ( 20.12, 67.80) {$0.05$};

\node[text=drawColor,anchor=base east,inner sep=0pt, outer sep=0pt, scale=  0.73] at ( 20.12, 82.43) {$0.4$};
\end{scope}
\begin{scope}
\path[clip] (  0.00,  0.00) rectangle (102.62, 93.95);
\definecolor{drawColor}{gray}{0.20}

\path[draw=drawColor,line width= 0.6pt,line join=round] ( 22.32, 37.89) --
	( 25.07, 37.89);

\path[draw=drawColor,line width= 0.6pt,line join=round] ( 22.32, 54.10) --
	( 25.07, 54.10);

\path[draw=drawColor,line width= 0.6pt,line join=round] ( 22.32, 70.30) --
	( 25.07, 70.30);

\path[draw=drawColor,line width= 0.6pt,line join=round] ( 22.32, 84.93) --
	( 25.07, 84.93);
\end{scope}
\begin{scope}
\path[clip] (  0.00,  0.00) rectangle (102.62, 93.95);
\definecolor{drawColor}{gray}{0.20}

\path[draw=drawColor,line width= 0.6pt,line join=round] ( 28.34, 24.92) --
	( 28.34, 27.67);

\path[draw=drawColor,line width= 0.6pt,line join=round] ( 50.18, 24.92) --
	( 50.18, 27.67);

\path[draw=drawColor,line width= 0.6pt,line join=round] ( 72.01, 24.92) --
	( 72.01, 27.67);

\path[draw=drawColor,line width= 0.6pt,line join=round] ( 93.85, 24.92) --
	( 93.85, 27.67);
\end{scope}
\begin{scope}
\path[clip] (  0.00,  0.00) rectangle (102.62, 93.95);
\definecolor{drawColor}{RGB}{0,0,0}

\node[text=drawColor,anchor=base,inner sep=0pt, outer sep=0pt, scale=  0.79] at ( 28.34, 17.27) {$1\%$};

\node[text=drawColor,anchor=base,inner sep=0pt, outer sep=0pt, scale=  0.79] at ( 50.18, 17.27) {$4\%$};

\node[text=drawColor,anchor=base,inner sep=0pt, outer sep=0pt, scale=  0.79] at ( 72.01, 17.27) {$7\%$};

\node[text=drawColor,anchor=base,inner sep=0pt, outer sep=0pt, scale=  0.79] at ( 93.85, 17.27) {$10\%$};
\end{scope}
\begin{scope}
\path[clip] (  0.00,  0.00) rectangle (102.62, 93.95);
\definecolor{drawColor}{RGB}{0,0,0}

\node[text=drawColor,anchor=base,inner sep=0pt, outer sep=0pt, scale=  0.85] at ( 61.10,  7.15) {Noise Level ($\sigma$)};
\end{scope}
\begin{scope}
\path[clip] (  0.00,  0.00) rectangle (102.62, 93.95);
\definecolor{drawColor}{RGB}{0,0,0}

\path[draw=drawColor,line width= 0.6pt,line join=round] ( 41.05, 41.30) -- ( 47.79, 41.30);
\end{scope}
\begin{scope}
\path[clip] (  0.00,  0.00) rectangle (102.62, 93.95);
\definecolor{drawColor}{RGB}{0,0,0}

\path[draw=drawColor,line width= 0.4pt,line join=round,line cap=round] ( 42.99, 39.87) rectangle ( 45.85, 42.72);
\end{scope}
\begin{scope}
\path[clip] (  0.00,  0.00) rectangle (102.62, 93.95);
\definecolor{drawColor}{RGB}{0,0,0}

\path[draw=drawColor,line width= 0.6pt,line join=round] ( 41.05, 41.30) -- ( 47.79, 41.30);
\end{scope}
\begin{scope}
\path[clip] (  0.00,  0.00) rectangle (102.62, 93.95);
\definecolor{drawColor}{RGB}{0,0,0}

\path[draw=drawColor,line width= 0.4pt,line join=round,line cap=round] ( 42.99, 39.87) rectangle ( 45.85, 42.72);
\end{scope}
\begin{scope}
\path[clip] (  0.00,  0.00) rectangle (102.62, 93.95);
\definecolor{drawColor}{RGB}{160,32,240}

\path[draw=drawColor,line width= 0.6pt,line join=round] ( 41.05, 32.86) -- ( 47.79, 32.86);
\end{scope}
\begin{scope}
\path[clip] (  0.00,  0.00) rectangle (102.62, 93.95);
\definecolor{drawColor}{RGB}{160,32,240}

\path[draw=drawColor,line width= 0.4pt,line join=round,line cap=round] ( 42.99, 31.44) rectangle ( 45.85, 34.29);
\end{scope}
\begin{scope}
\path[clip] (  0.00,  0.00) rectangle (102.62, 93.95);
\definecolor{drawColor}{RGB}{160,32,240}

\path[draw=drawColor,line width= 0.6pt,line join=round] ( 41.05, 32.86) -- ( 47.79, 32.86);
\end{scope}
\begin{scope}
\path[clip] (  0.00,  0.00) rectangle (102.62, 93.95);
\definecolor{drawColor}{RGB}{160,32,240}

\path[draw=drawColor,line width= 0.4pt,line join=round,line cap=round] ( 42.99, 31.44) rectangle ( 45.85, 34.29);
\end{scope}
\begin{scope}
\path[clip] (  0.00,  0.00) rectangle (102.62, 93.95);
\definecolor{drawColor}{RGB}{0,0,0}

\node[text=drawColor,anchor=base west,inner sep=0pt, outer sep=0pt, scale=  0.88] at ( 50.80, 38.27) {$|d_1|$};
\end{scope}
\begin{scope}
\path[clip] (  0.00,  0.00) rectangle (102.62, 93.95);
\definecolor{drawColor}{RGB}{0,0,0}

\node[text=drawColor,anchor=base west,inner sep=0pt, outer sep=0pt, scale=  0.6] at ( 50.80, 29.83) {$\trsp{\hat{\bw}_0} \bQ_1 \hat{\bw}_0$+$\| \bQ_1 \hat{\bw}\|_2$};
\end{scope}
\end{tikzpicture} \label{fig:di} }
	\caption{Numerical illustration of condition \eqref{eq:noisy+outlier-free2-condition} of \cref{theorem:noisy+outlier-free2} ($500$ trials, $\ell=100$).}
\end{figure}
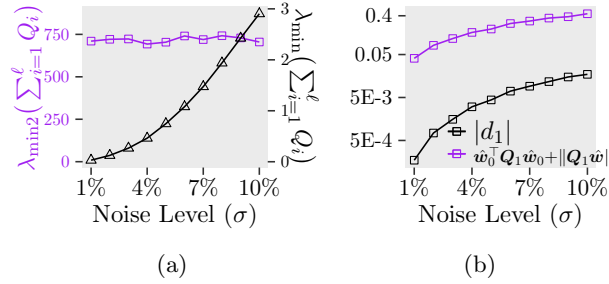

\cref{theorem:noisy+outlier-free2} is better understood via numerics. We take randomly generated $\ell=100$ point pairs $(\by_i,\bx_i)$'s with $\bx_i\sim \cN(0,\bI_3)$, and add different levels of Gaussian noise $\epsilon_i\sim \cN(0,\sigma^2\bI_3)$, where $\sigma$ ranges from $1\%$ to $10\%$. The values of $\lambda_{\textnormal{min}}\big( \sum_{i=1}^{\ell} \bQ_i \big)$ and $\lambda_{\textnormal{min2}}\big(  \sum_{i=1}^{\ell} \bQ_i \big)$, and thus $\zeta$, are shown in Figure \ref{fig:eigengap}, where one might observe that $\zeta\approx 250$ for $10\%$ noise, $\zeta\approx 25000$ for $1\%$ noise, and, in general, $\zeta=\infty$ for the noiseless case. This empirically validates the assumption $\zeta\geq \ell/(\ell-1)=100/99$.

We then elaborate the more complicated condition \eqref{eq:noisy+outlier-free2-condition}. First we recall that $c_i^2 > \trsp{\hat{\bw}_0} \bQ_i \hat{\bw}_0$ is essential for \eqref{eq:TLS-Q} to preserve all inliers. Second, we argue that the term $\| \bQ_i \hat{\bw}\|_2$ in \eqref{eq:noisy+outlier-free2-condition} is also essential, as it accounts for the fact that noise destroys the inequality $\lambda_{\textnormal{min}}(\bQ_i)  - \trsp{\hat{\bw}_0} \bQ_i \hat{\bw}_0\geq 0$ which holds in the noiseless case (where $\lambda_{\textnormal{min}}(\bQ_i)  = \trsp{\hat{\bw}_0} \bQ_i \hat{\bw}_0=0$) but gets violated (in general) in the presence of noise. Finally, \eqref{eq:noisy+outlier-free2-condition} also incurs a curious term $(|d_i|+d_i)/2$, with $d_i$ defined in a sophisticated way \eqref{eq:noisy+outlier-free2-condition}. If $d_i<0$ then this term is $0$. Thus it remains to understand the values of $|d_i|$. In particular, we plotted the values of $|d_1|$ in Figure \ref{fig:di} in comparison to $\trsp{\hat{\bw}_0} \bQ_1 \hat{\bw}_0 + \| \bQ_1 \hat{\bw}\|_2$, and observed that $|d_1|$ is two orders of magnitude smaller (there is nothing special about the choice of index $1$). In fact, as noise approaches zero, we have $\bQ_i\hat{\bw}_0\to 0$ and (in general) $\zeta\to \infty$, hence $d_i\to 0$ by definition \eqref{eq:noisy+outlier-free2-condition}. Overall, condition \eqref{eq:noisy+outlier-free2-condition} degenerates into $c_i^2>0$ in the noiseless case. 

While the term $(|d_i|+d_i)/2$ is quite small (Figure \ref{fig:di}) and sometimes harmless (e.g., when $d_i<0$), it appears as an artifact of our analysis, and we expect an ideal condition for the noisy + outlier-free case to be $c_i^2 > \trsp{\hat{\bw}_0} \bQ_i \hat{\bw}_0 + \norm{\bQ_i \hat{\bw}_0}$. However, proof under this alternative condition demands showing some matrix inequality that involves a sum of matrix inverses always holds; we were not able to prove it.

Finally, we discuss the error bound \eqref{eq:dist-eigvector2gt}. It becomes zero as $\epsilon_i\to 0$ and thus $\hat{\bw}_0=\bw_0^*$, provided that $\|\bx_i \|_2 / \lambda_{\textnormal{min2}}\big(\sum_{i=1}^{\ell} \bP_i\big)$ is not too large. The denominator $\lambda_{\textnormal{min2}}\big(\sum_{i=1}^{\ell} \bP_i\big)$ seems inevitable, as it usually determines the stability of solving a minimum eigenvalue problem: If $\lambda_{\textnormal{min2}}\big(\sum_{i=1}^{\ell} \bP_i\big) \to 0=\lambda_{\textnormal{min}}\big(\sum_{i=1}^{\ell} \bP_i\big)$, then $\hat{\bw}_0$ can be arbitrarily far from $\bw_0^*$ even in the slightest presence of noise. That said, $\|\bx_i \|_2 / \lambda_{\textnormal{min2}}\big(\sum_{i=1}^{\ell} \bP_i\big)$ is actually quite small and well-behaved, at least for random Gaussian data; see \cref{section:probability-bounds} for details.

\subsection{The Noisy + Outliers Case}\label{subsection:noisy+outliers}
With the proof ideas of Theorems \ref{theorem:noiseless+outlier+c} and \ref{theorem:noisy+outlier-free2}, we get:
\begin{theorem}[Noisy Point Sets with Outliers]\label{theorem:noisy+outliers}
	Let $\bQ_1,\dots,\bQ_{k^*}$ be inliers, the rest $\bQ_j$'s outliers, and $\hat{\bw}_0$ a global minimizer of  \eqref{eq:TLS-Q}. Define 
	\begin{equation}\label{eq:eigengap2}
		\zeta_{\textnormal{in}} := \frac{\lambda_{\textnormal{min2}}\big(  \sum_{i=1}^{k^*} 	\bQ_i \big)}{\lambda_{\textnormal{min}}\big(  \sum_{i=1}^{k^*} 	\bQ_i \big)}.
	\end{equation}			
	Assume (1) $\zeta_{\textnormal{in}} \geq k^*/(k^*-1)$, (2) for every $j=k^*+1,\dots,\ell$, we have $0<c_j^2<\lambda_{\textnormal{min}}(\bQ_{j})$, (3) for every $i=1,\dots,k^*$, \eqref{eq:noisy+outlier-free2-condition} holds with $d_i$ now defined as
	\begin{equation}
		d_i := \frac{\sum_{i=1}^{k^*} \trsp{\hat{\bw}_0}\bQ_i\hat{\bw}_0}{k^*} - \trsp{\hat{\bw}_0} \bQ_i \hat{\bw}_0 + \frac{\lambda_{\textnormal{max}}\Big(\sum_{j\neq i} \big(  \bQ_{i}  -  \bQ_{j} \big)\Big)}{\zeta(k^* - 1)}. \label{eq:di2}
	\end{equation}
	Then \eqref{eq:SDP} is tight and, similarly to \eqref{eq:dist-eigvector2gt} we have
	\begin{equation}\label{eq:dist-eigvector2gt2}
		\sin^2(\hat{\tau}^*_0) \leq \frac{4\sum_{i=1}^{k^*}\norm{\beps_i} \norm{\bx_i}}{\lambda_{\textnormal{min2}}\big(\sum_{i=1}^{k^*} \bP_i\big)}, \ \ \ \ \  \sin^2(\hat{\tau}^*_0):= 1 - (\trsp{\hat{\bw}_0}\bw_0^*)^2
	\end{equation}
	Here we recall that $\bw_0^*\in\bbS^3$ is the ground-truth unit quaternion, and each $\bP_i$ is the ``pure data'' part of $\bQ_i$ that satisfies $\bP_i\bw_0^*=0,\bP_i\succeq0$ (\cref{lemma:rotation-quaternion2}). 
\end{theorem}
\begin{proof}
	The given assumptions ensure that \eqref{eq:TLS-Q} rejects all outliers and admit all inliers, and the minimum of \eqref{eq:TLS-Q} is $\sum_{i=1}^{k^*} \trsp{\hat{\bw}_0}\bQ_i \hat{\bw}_0  + \sum_{j=k^*+1}^{\ell} c_j^2$.  For $i=1,\dots,k^*$ let $\pD_{0i}$ be defined as in \eqref{eq:D0i-Si}, and for $j=k^*+1,\dots,\ell$ let $[\hat{\pD}]_{0j}:= \frac{1}{2}(\bQ_j - c_j^2\bI_4)$. Let $\hat{\mu} := \sum_{i=1}^{k^*} \big(\trsp{\hat{\bw}_0}\bQ_i \hat{\bw}_0 - c_i^2\big)$. One then verifies the optimality conditions \eqref{eq:stationary2} and \eqref{eq:value-mu2} of \cref{prop:tightnesscondition} are satisfied. \eqref{eq:dualfeasibility2} is equivalent to ($\forall \bz_i\in\bbR^4$)
	\begin{equation*}
		\underbrace{-\hat{\mu} \norm{\bz_0}^2 + \sum_{i=1}^{k^*} \Big( 2\trsp{\bz_i}[\hat{\pD}]_{0i} \bz_i -\trsp{\bz_0}\big(2[\hat{\pD}]_{0i} - \bQ_i+c_i^2\bI_4\big)\bz_i \Big) }_{\textnormal{Inlier Term}} + \underbrace{\sum_{i=k^*+1}^{\ell}\bz_j(\bQ_j - c_j^2\bI_4)\bz_j}_{\textnormal{Outlier Term}} \geq 0,
	\end{equation*} 
	Since $\bQ_j\succeq c_j^2\bI_4$, the outlier term is non-negative. Under the given assumptions, one can replace $\ell$ by $k^*$ in the proof of \cref{theorem:noisy+outlier-free2} and then find the inlier term is also non-negative. This finishes proving \eqref{eq:dualfeasibility2} and thus the tightness of \eqref{eq:SDP}. The error bound \eqref{eq:dist-eigvector2gt2} follows from the proof of \cref{theorem:noisy+outlier-free2} with $\ell$ also replaced by $k^*$.
\end{proof}

\section{Probabilistic Interpretation of Error Bound \eqref{eq:dist-eigvector2gt}}\label{section:probability-bounds}
\begin{figure}[t]
	\centering
\begin{tikzpicture}[x=1pt,y=1pt]
\definecolor{fillColor}{RGB}{255,255,255}
\path[use as bounding box,fill=fillColor,fill opacity=0.00] (0,0) rectangle (144.54, 93.95);
\begin{scope}
\path[clip] (  0.00,  0.00) rectangle (144.54, 93.95);
\definecolor{drawColor}{RGB}{255,255,255}
\definecolor{fillColor}{RGB}{255,255,255}

\path[draw=drawColor,line width= 0.6pt,line join=round,line cap=round,fill=fillColor] (  0.00,  0.00) rectangle (144.54, 93.95);
\end{scope}
\begin{scope}
\path[clip] ( 21.34, 29.01) rectangle (139.04, 88.45);
\definecolor{fillColor}{gray}{0.92}

\path[fill=fillColor] ( 21.34, 29.01) rectangle (139.04, 88.45);
\definecolor{drawColor}{RGB}{255,0,0}

\path[draw=drawColor,line width= 0.6pt,line join=round] ( 26.69, 68.84) --
	( 44.52, 50.38) --
	( 62.36, 46.65) --
	( 80.19, 45.66) --
	( 98.02, 45.36) --
	(115.86, 45.27) --
	(133.69, 45.24);

\path[draw=drawColor,line width= 0.4pt,line join=round,line cap=round] ( 25.26, 67.41) rectangle ( 28.12, 70.27);

\path[draw=drawColor,line width= 0.4pt,line join=round,line cap=round] ( 43.10, 48.96) rectangle ( 45.95, 51.81);

\path[draw=drawColor,line width= 0.4pt,line join=round,line cap=round] ( 60.93, 45.22) rectangle ( 63.78, 48.08);

\path[draw=drawColor,line width= 0.4pt,line join=round,line cap=round] ( 78.76, 44.24) rectangle ( 81.62, 47.09);

\path[draw=drawColor,line width= 0.4pt,line join=round,line cap=round] ( 96.60, 43.94) rectangle ( 99.45, 46.79);

\path[draw=drawColor,line width= 0.4pt,line join=round,line cap=round] (114.43, 43.84) rectangle (117.28, 46.69);

\path[draw=drawColor,line width= 0.4pt,line join=round,line cap=round] (132.26, 43.81) rectangle (135.12, 46.66);
\definecolor{drawColor}{RGB}{0,0,0}

\path[draw=drawColor,line width= 0.6pt,dash pattern=on 4pt off 4pt ,line join=round] ( 26.69, 31.72) --
	( 44.52, 31.72) --
	( 62.36, 31.72) --
	( 80.19, 31.72) --
	( 98.02, 31.72) --
	(115.86, 31.72) --
	(133.69, 31.72);

\path[draw=drawColor,line width= 0.6pt,dash pattern=on 4pt off 4pt ,line join=round] ( 26.69, 85.75) --
	( 44.52, 85.75) --
	( 62.36, 85.75) --
	( 80.19, 85.75) --
	( 98.02, 85.75) --
	(115.86, 85.75) --
	(133.69, 85.75);

\path[draw=drawColor,line width= 0.6pt,dash pattern=on 4pt off 4pt ,line join=round] ( 26.69, 45.22) --
	( 44.52, 45.22) --
	( 62.36, 45.22) --
	( 80.19, 45.22) --
	( 98.02, 45.22) --
	(115.86, 45.22) --
	(133.69, 45.22);
\end{scope}
\begin{scope}
\path[clip] (  0.00,  0.00) rectangle (144.54, 93.95);
\definecolor{drawColor}{RGB}{0,0,0}

\node[text=drawColor,anchor=base east,inner sep=0pt, outer sep=0pt, scale=  0.73] at ( 16.39, 29.22) {$1/4$};

\node[text=drawColor,anchor=base east,inner sep=0pt, outer sep=0pt, scale=  0.73] at ( 16.39, 42.72) {$3/8$};

\node[text=drawColor,anchor=base east,inner sep=0pt, outer sep=0pt, scale=  0.73] at ( 16.39, 83.25) {$3/4$};
\end{scope}
\begin{scope}
\path[clip] (  0.00,  0.00) rectangle (144.54, 93.95);
\definecolor{drawColor}{gray}{0.20}

\path[draw=drawColor,line width= 0.6pt,line join=round] ( 18.59, 31.72) --
	( 21.34, 31.72);

\path[draw=drawColor,line width= 0.6pt,line join=round] ( 18.59, 45.22) --
	( 21.34, 45.22);

\path[draw=drawColor,line width= 0.6pt,line join=round] ( 18.59, 85.75) --
	( 21.34, 85.75);
\end{scope}
\begin{scope}
\path[clip] (  0.00,  0.00) rectangle (144.54, 93.95);
\definecolor{drawColor}{gray}{0.20}

\path[draw=drawColor,line width= 0.6pt,line join=round] ( 26.69, 26.26) --
	( 26.69, 29.01);

\path[draw=drawColor,line width= 0.6pt,line join=round] ( 44.52, 26.26) --
	( 44.52, 29.01);

\path[draw=drawColor,line width= 0.6pt,line join=round] ( 62.36, 26.26) --
	( 62.36, 29.01);

\path[draw=drawColor,line width= 0.6pt,line join=round] ( 80.19, 26.26) --
	( 80.19, 29.01);

\path[draw=drawColor,line width= 0.6pt,line join=round] ( 98.02, 26.26) --
	( 98.02, 29.01);

\path[draw=drawColor,line width= 0.6pt,line join=round] (115.86, 26.26) --
	(115.86, 29.01);

\path[draw=drawColor,line width= 0.6pt,line join=round] (133.69, 26.26) --
	(133.69, 29.01);
\end{scope}
\begin{scope}
\path[clip] (  0.00,  0.00) rectangle (144.54, 93.95);
\definecolor{drawColor}{RGB}{0,0,0}

\node[text=drawColor,anchor=base west,inner sep=0pt, outer sep=0pt, scale=  0.79] at ( 21.35, 17.27) {10};

\node[text=drawColor,anchor=base west,inner sep=0pt, outer sep=0pt, scale=  0.55] at ( 29.26, 20.51) {1};

\node[text=drawColor,anchor=base west,inner sep=0pt, outer sep=0pt, scale=  0.79] at ( 39.18, 17.27) {10};

\node[text=drawColor,anchor=base west,inner sep=0pt, outer sep=0pt, scale=  0.55] at ( 47.10, 20.51) {2};

\node[text=drawColor,anchor=base west,inner sep=0pt, outer sep=0pt, scale=  0.79] at ( 57.01, 17.27) {10};

\node[text=drawColor,anchor=base west,inner sep=0pt, outer sep=0pt, scale=  0.55] at ( 64.93, 20.51) {3};

\node[text=drawColor,anchor=base west,inner sep=0pt, outer sep=0pt, scale=  0.79] at ( 74.85, 17.27) {10};

\node[text=drawColor,anchor=base west,inner sep=0pt, outer sep=0pt, scale=  0.55] at ( 82.76, 20.51) {4};

\node[text=drawColor,anchor=base west,inner sep=0pt, outer sep=0pt, scale=  0.79] at ( 92.68, 17.27) {10};

\node[text=drawColor,anchor=base west,inner sep=0pt, outer sep=0pt, scale=  0.55] at (100.60, 20.51) {5};

\node[text=drawColor,anchor=base west,inner sep=0pt, outer sep=0pt, scale=  0.79] at (110.51, 17.27) {10};

\node[text=drawColor,anchor=base west,inner sep=0pt, outer sep=0pt, scale=  0.55] at (118.43, 20.51) {6};

\node[text=drawColor,anchor=base west,inner sep=0pt, outer sep=0pt, scale=  0.79] at (128.35, 17.27) {10};

\node[text=drawColor,anchor=base west,inner sep=0pt, outer sep=0pt, scale=  0.55] at (136.26, 20.51) {7};
\end{scope}
\begin{scope}
\path[clip] (  0.00,  0.00) rectangle (144.54, 93.95);
\definecolor{drawColor}{RGB}{0,0,0}

\node[text=drawColor,anchor=base,inner sep=0pt, outer sep=0pt, scale=  0.85] at ( 80.19,  7.15) {$\ell$};
\end{scope}
\end{tikzpicture} 
	\caption{Numerical values of $\frac{\sum_{i=1}^{\ell} \| \bx_i\|_2^2}{\lambda_{\textnormal{min2}} (\sum_{i=1}^{\ell}\bP_i )}$ with $\bx_i\sim \cN(0, \bI_3)$ in red, $100$ trials. \label{fig:ratio}}
\end{figure}

Since the entries of $\bP_i$ depend quadratically on entries of $\bx_i$,  we naturally expect  $\lambda_{\textnormal{min2}} \big(\sum_{i=1}^{\ell}\bP_i \big)=\Theta(\sum_{i=1}^{\ell} \| \bx_i \|_2^2)$ for random points $\bx_i$'s; \cref{prop:error_bound} confirms that this is indeed true with high probability:

\begin{proposition}\label{prop:error_bound}
	Suppose that every 3D point $\bx_i$ has $\cN(0,1)$ entries. Let $t$ be some positive constant and let $\ell$ be large enough in the sense that
	\begin{equation}\label{eq:condition-ell}
		\ell \geq  \big(4 + 2\sqrt{3} + 2t\big) \sqrt{\ell} +  \big(\sqrt{3} + t\big)^2.
	\end{equation}
	Then the following holds with probability at least $1-\exp(-t^2/2) - 2\exp(-3t^2/8)$:
	\begin{equation} \label{eq:x^2=minlambda2}
		\frac{1}{4}  \leq \frac{\sum_{i=1}^{\ell} \norm{\bx_i}^2}{\lambda_{\textnormal{min2}} \big(\sum_{i=1}^{\ell}\bP_i \big)} \leq   \frac{3}{4}
	\end{equation}
	Recall that $\bP_i $ is the pure data part of $\bQ_i$ in the noisy case (cf. \cref{lemma:rotation-quaternion2}). 
\end{proposition}
\begin{remark}[Asymptotic Behavior]
	With \cref{lemma:Wainwright-2.11,lemma:Wainwright-6.1,lemma:lambdamin=lambdamax} and \cite{Wainwright-book2019}, we get
	\begin{equation}
		\ell \to\infty \ \ \ \Rightarrow  \ \ \ \frac{\sum_{i=1}^{\ell} \norm{\bx_i}^2}{\lambda_{\textnormal{min2}} \big(\sum_{i=1}^{\ell}\bP_i \big)} = \frac{3}{8}
	\end{equation}
\end{remark}
\begin{proof}[Proof of \cref{prop:error_bound}]
	Since $\bx_i$'s has i.i.d.  $\cN(0,1)$ entries, $\lambda_{\textnormal{min2}} \Big(\sum_{i=1}^{\ell}\bP_i \Big)$ is not equal to zero with probability $1$. Then \cref{lemma:lambdamin=lambdamax} implies that the left bound holds. For the right bound, applying the union bound with \cref{lemma:Wainwright-2.11,lemma:Wainwright-6.1,lemma:lambdamin=lambdamax}, it holds that
	\begin{equation}\label{eq:x^2=minlambda2_ub}
		\frac{\sum_{i=1}^{\ell} \norm{\bx_i}^2}{\lambda_{\textnormal{min2}} \Big(\sum_{i=1}^{\ell}\bP_i \Big)} \leq   \frac{3\ell + 3\sqrt{\ell}}{4\cdot (2\ell - 3\sqrt{\ell} -  2\big(\sqrt{3} + t\big) \sqrt{\ell} -  \big(\sqrt{3} + t\big)^2 )}
	\end{equation}
	with probability  at least $1-\exp(-t^2/2) - 2\exp(-3t^2/8)$. Next, from \eqref{eq:condition-ell} we get
	\begin{equation}
		2\ell - 3\sqrt{\ell} -  2\big(\sqrt{3} + t\big) \sqrt{\ell} -  \big(\sqrt{3} + t\big)^2 \geq \ell + \sqrt{\ell},
	\end{equation}
	which means the right-hand side of \eqref{eq:x^2=minlambda2_ub} does not exceed $3/4$. The proof is finished.
\end{proof}

\begin{lemma}[Example 2.11 of \cite{Wainwright-book2019}]\label{lemma:Wainwright-2.11}
	If each $\bx_i\in\bbR^3$ has i.i.d. $\cN(0,1)$ entries, then it holds for some $t > 0$ with probability at least $1-2\exp(-3 t^2/8)$ that
	\begin{equation*}
		3\ell -3 \sqrt{\ell} \leq \sum_{i=1}^{\ell} \norm{\bx_i}^2 \leq 3\ell+ 3 \sqrt{\ell}, \ \ \ \ \ \ t > 0
	\end{equation*}
\end{lemma}

\begin{lemma}[Theorem 6.1 and Example 6.2 of \cite{Wainwright-book2019}]\label{lemma:Wainwright-6.1} 
	If each $\bx_i\in\bbR^3$ has i.i.d. $\cN(0,1)$ entries, then for some $t > 0$ with probability at least $1-\exp(- t^2/2)$ we have
	\begin{equation*}
		\lambda_{\textnormal{max}}\Big( \sum_{i=1}^{\ell} \bx_i \trsp{\bx_i} \Big) \leq \ell + 2\big(\sqrt{3} + t\big) \sqrt{\ell} +  \big(\sqrt{3} + t\big)^2 
	\end{equation*}
\end{lemma}

\begin{lemma}\label{lemma:lambdamin=lambdamax}
	Recall that $\bP_i$ is the pure data part of $\bQ_i$ in the noisy case (cf. \cref{lemma:rotation-quaternion2}), and that $\lambda_{\textnormal{min2}}(\cdot)$ denotes the second smallest eigenvalue of a matrix. We have
	\begin{equation}\label{eq:lambdamin=lambdamax}
		\lambda_{\textnormal{min2}} \Big(\sum_{i=1}^{\ell}\bP_i \Big)= \sum_{i=1}^{\ell} 4\norm{\bx_i}^2 - 
		\lambda_{\textnormal{max}} \Big( \sum_{i=1}^{\ell} 4\bx_i \trsp{\bx_i} \Big) 
	\end{equation}
\end{lemma}
\begin{proof}[Proof of \cref{lemma:lambdamin=lambdamax}]
	With any $\bw_2\in\bbS^3$, $\phi$ and $\bm{b}$ the rotation angle and rotation axis of $\trsp{\bR_2}\bR^*_0$ respectively, and with \cref{lemma:rotation-quaternion3}, we have the following equivalence:
	\begin{equation*}
		\trsp{\bw_2} \bw^*_0=0  \Leftrightarrow \phi = \pi  \Leftrightarrow  \trsp{\bR_2}\bR^*_0 = 2\bm{b}\trsp{\bm{b}} - \bI_3
	\end{equation*}
	In the last equivalence we used the axis-angle representation of 3D rotations
	\begin{equation*}
		\trsp{\bR_2}\bR^*_0 = \bm{b}\trsp{\bm{b}} + [\bm{b}]_{\times} \sin(\phi) + \big( \bI_3 - \bm{b} \trsp{\bm{b}} \big) \cos(\phi)
	\end{equation*}
	and substituted $\phi=\pi$ into it; here, $[\bm{b}]_{\times}\in\bbR^{3\times 3}$ is the cross product matrix that satisfies $[\bm{b}]_{\times}\bq = \bm{b}\times \bq$ for any $\bq\in\bbR^3$. On the other hand, \cref{lemma:rotation-quaternion} implies
	\begin{equation*}
		\trsp{\bw_2} \Big( \sum_{i=1}^{\ell}\bP_i  \Big) \bw_2= \sum_{i=1}^{\ell}\norm{\bR_0^*\bx_i - \bR_2 \bx_i}^2 = \sum_{i=1}^{\ell} 2\norm{\bx_i}^2 - \sum_{i=1}^{\ell} 2\trsp{\bx_i} \trsp{\bR_2}\bR_0^* \bx_i.
	\end{equation*}
	Clearly, as the unit quaternion $\bw_2$ varies in $\bbS^3$ with $\trsp{\bw_2}\bw^*_0=0$, the axis of $\bm{b}$ of $\trsp{\bR_2}\bR^*_0$ can be any element of $\bbS^2$, and vice versa. Combine the above, and we arrive at
	\begin{equation}\label{eq:lambdamin=lambdamax2}
		\begin{split}
			\min_{\substack{\bw\in\bbS^3 \\ \trsp{\bw_2}\bw^*_0=0}} \trsp{\bw_2} \Big( \sum_{i=1}^{\ell}\bP_i  \Big) \bw_2 &= \sum_{i=1}^{\ell} 2\norm{\bx_i}^2 - \max_{\bm{b}\in\bbS^2} \sum_{i=1}^{\ell} 2\trsp{\bx_i} \big( 2\bm{b}\trsp{\bm{b}} - \bI_3 \big) \bx_i 
		\end{split}
	\end{equation}
	Since $\bw_0^*$ corresponds to the minimum eigenvalue $0$ of $ \sum_{i=1}^{\ell}\bP_i $, we know that the left-hand side of \eqref{eq:lambdamin=lambdamax2} is equal to $\lambda_{\textnormal{min2}} \big( \sum_{i=1}^{\ell}\bP_i  \big)$. Also observe that the right-hand side of \eqref{eq:lambdamin=lambdamax2} is equal to that of \eqref{eq:lambdamin=lambdamax}, and we have finished the proof.
\end{proof}

\begin{lemma}\label{lemma:rotation-quaternion3}
	For two 3D rotations $\bR_1$ and $\bR_2$, their repsective unit quaternion representations $\bw_1$ and $\bw_2$, and the rotation angle $\phi_{12}\in[0,2\pi]$ of $\trsp{\bR_1}\bR_2$, we have 
	\begin{equation*}
		\trsp{\bw_1}\bw_2 = \cos\Big(\frac{\phi_{12}}{2}\Big).
	\end{equation*}
	In particular, we get that $\trsp{\bw_1}\bw_2 = 0$ if and only if $\phi_{12}=\pi$.
\end{lemma}
\begin{proof}[Proof of \cref{lemma:rotation-quaternion3}]
	Denote by $\bm{b}_1$ (resp. $\bm{b}_2$) and $\phi_1$ (resp. $\phi_2$) the rotation axis and angle of $\bR_1$ (resp. $\bR_2$). Then, from a well-known relation between unit quaternions and axis-angle representation of rotations, we know that $\bw_1 = [\cos (\frac{\phi_1}{2});  \bm{b}_1 \sin (\frac{\phi_1}{2})]$ and  $\bw_2 = [\cos (\frac{\phi_2}{2});  \bm{b}_2 \sin (\frac{\phi_2}{2})]$, and therefore we have
	\begin{equation*}
		\trsp{\bw_1}\bw_2 = \cos \Big(\frac{\phi_1}{2} \Big)\cdot \cos \Big(\frac{\phi_2}{2} \Big) +  \trsp{\bm{b}_1}\bm{b}_2 \cdot \sin \Big(\frac{\phi_1}{2} \Big) \cdot \sin \Big(\frac{\phi_2}{2} \Big).
	\end{equation*}
	Since $\trsp{\bR_1}$ is a 3D rotation with angle $\phi_1$ and axis $-\bm{b}_1$, it is well known that the composition $\trsp{\bR_1}\bR_2$ of two 3D rotations has its rotation angle $\phi_{12}$ satisfying
	\begin{equation}\label{eq:rotation-composition}
		\cos\Big(\frac{\phi_{12}}{2}\Big) = \cos \Big(\frac{\phi_1}{2} \Big)\cdot \cos \Big(\frac{\phi_2}{2} \Big) - \trsp{(-\bm{b}_1)}\bm{b}_2 \cdot \sin \Big(\frac{\phi_1}{2}\Big) \cdot \sin \Big(\frac{\phi_2}{2}\Big).
	\end{equation}
	Equality \eqref{eq:rotation-composition} can be proved via spherical trigonometry, or \eqref{eq:rotation-quaternion}, or Rodrigues' rotation formula, with some algebraic manipulations. We have thus finished the proof.
\end{proof}

\section{Discussion and Future Work}\label{section:discussion}
We have investigated the tightness of a semidefinite relaxation \eqref{eq:SDP} of truncated least-squares for robust rotation search in four different cases, and in each case we either showed improvements over prior work or proved new theoretical results. Our investigation can potentially be borrowed to understand semidefinite relaxations of many other geometric vision tasks; see \cite{Yang-arXiv2021b} for $6$ examples of truncated least-squares and see also \cite{Cifuentes-MP2021,Alfassi-Sensors2021,Brynte-JMIV2022}. 

As is common in the optimization literature, the relaxation we analyzed is at the first (i.e., lowest) relaxation order of the \textit{Lasserre} hierarchy \cite{Lasserre-SIAM-J-O2001}, or otherwise known as the \textit{Shor relaxation} \cite{Shor-AOR1990}. A tighter relaxation that has quadratically more constraints than \eqref{eq:SDP} exists (cf. \cite{Yang-ICCV2019}). However, analyzing this tighter relaxation is significantly harder, as one needs to either (1) construct quadratically more dual certificates during the proof, or (2) use more abstract optimality conditions (cf. \cite{Cifuentes-SIAM-MAA2021,Cifuentes-MP2021}). Therefore, we leave this challenging question to future work.

\appendix

\section{Proof of \cref{lemma:property-Si}}
\begin{proof}[Proof of \cref{lemma:property-Si}]
	We first prove the equality $\sum_{i=1}^{\ell} \widehat{\bS}_i = \sum_{i=1}^{\ell} \big( \bQ_i  -  \trsp{\hat{\bw}_0}\bQ_i\hat{\bw}_0 \bI_4  \big)$. First, by definition of $\widehat{\bT}_i$ \eqref{eq:Ti} and simple calculation, we have
	\begin{equation*}
		\sum_{i=1}^{\ell} \widehat{\bT}_i = \sum_{i=1}^{\ell} \big( \trsp{\widehat{\bV}_0}\bQ_i \widehat{\bV}_0 -  \trsp{\hat{\bw}_0}\bQ_i\hat{\bw}_0 \bI_3  \big)
	\end{equation*}
	Then, with the definition of $\widehat{\bS}$ \eqref{eq:D0i-Si}, we arrive at the following equivalence:
	\begin{equation*}
		\begin{split}
			&\sum_{i=1}^{\ell} \widehat{\bS}_i = \sum_{i=1}^{\ell} \big( \bQ_i  -  \trsp{\hat{\bw}_0}\bQ_i\hat{\bw}_0 \bI_4  \big)  \\
			\Leftrightarrow & \sum_{i=1}^{\ell} \begin{bmatrix}
				\widehat{\bT}_i  & 0 \\
				0 & 0
			\end{bmatrix} = \sum_{i=1}^{\ell} \big( \trsp{\widehat{\bV}}\bQ_i \widehat{\bV} -  \trsp{\hat{\bw}_0}\bQ_i\hat{\bw}_0 \bI_4  \big) \\
			\Leftrightarrow &  \begin{bmatrix}
				\sum_{i=1}^{\ell} \big( \trsp{\widehat{\bV}_0}\bQ_i \widehat{\bV}_0 -  \trsp{\hat{\bw}_0}\bQ_i\hat{\bw}_0 \bI_3  \big)  & 0 \\
				0 & 0
			\end{bmatrix} = \sum_{i=1}^{\ell} \big( \trsp{\widehat{\bV}}\bQ_i \widehat{\bV} -  \trsp{\hat{\bw}_0}\bQ_i\hat{\bw}_0 \bI_4  \big) \\
			\Leftrightarrow &  \begin{bmatrix}
				\sum_{i=1}^{\ell} \trsp{\widehat{\bV}_0}\bQ_i \widehat{\bV}_0   & 0 \\
				0 & 0
			\end{bmatrix} =  \begin{bmatrix}
				\sum_{i=1}^{\ell} \trsp{\widehat{\bV}_0}\bQ_i \widehat{\bV}_0   & & \sum_{i=1}^{\ell} \trsp{\widehat{\bV}_0}\bQ_i \hat{\bw}_0 \\
				\sum_{i=1}^{\ell} \trsp{\hat{\bw}_0}\bQ_i \widehat{\bV}_0 && 0
			\end{bmatrix}  
		\end{split}
	\end{equation*}
	Since $\trsp{\widehat{\bV}_0} \hat{\bw}_0=0$ and $\hat{\bw}_0$ is an eigenvector of $\sum_{i=1}^{\ell} \bQ_i$, the above equality holds.
	
	We now prove $\widehat{\bS}_{i}\succeq 0$. By the definition of  $\widehat{\bS}$ \eqref{eq:D0i-Si}, it suffices to show $\widehat{\bT}_{i}\succeq 0$. Now, using the definition of $\zeta$ \eqref{eq:eigengap} and the assumption $\zeta\geq \ell/(\ell-1)$, we have
	\begin{equation*}
		\begin{split}
			&\ \ \widehat{\bS}_{i}\succeq 0 \Leftrightarrow \widehat{\bT}_{i}\succeq 0\\
			\Leftrightarrow&\ \ \frac{\zeta-\frac{\ell}{\ell-1}}{\zeta}\trsp{\widehat{\bV}_0} \bQ_{i} \widehat{\bV}_0 + \frac{\sum_{j= 1}^\ell \trsp{\widehat{\bV}_0} \bQ_{j} \widehat{\bV}_0}{\zeta(\ell - 1)} - \frac{\Big( \sum_{i=1}^{\ell} \trsp{\hat{\bw}_0}\bQ_i\hat{\bw}_0 \Big) \bI_3}{\ell} \succeq 0 \\
			\Leftarrow& \ \ \frac{\sum_{j= 1}^\ell \trsp{\widehat{\bV}_0} \bQ_{j} \widehat{\bV}_0}{\zeta(\ell - 1)} - \frac{\lambda_{\textnormal{min}}\Big( \sum_{i=1}^{\ell} \bQ_i \Big) \bI_3}{\ell} \succeq 0 \\
			\Leftrightarrow &\ \ \sum_{j= 1}^\ell \trsp{\widehat{\bV}_0} \bQ_{j} \widehat{\bV}_0 -   \frac{\ell-1}{\ell}\lambda_{\textnormal{min2}}\Big( \sum_{i=1}^{\ell} \bQ_i \Big) \bI_3 \succeq 0,
		\end{split}
	\end{equation*}
	which holds true, as the minimum eigenvalue of $\sum_{j= 1}^\ell \trsp{\widehat{\bV}_0} \bQ_{j} \widehat{\bV}_0$ is exactly the second smallest eigenvalue of $\sum_{i=1}^{\ell} \bQ_i$ (by the definition of $\widehat{\bV}_0$). This proves $\widehat{\bS}_i\succeq 0$.
	
	Finally, it remains to prove the last inequality:
	\begin{equation}
		\begin{split}
			\widehat{\bS}_{i} + c_i^2\bI_4 -  \bQ_i \succ 0 \Leftrightarrow&\  \trsp{\bz_i} \Big(\widehat{\bS}_{i} + c_i^2\bI_4 -  \bQ_i\Big) \bz_i \geq 0,\  \forall \bz_i \in\bbR^4 \\
			\Leftrightarrow &\ \trsp{\ba_i} \bigg( \begin{bmatrix}
				\widehat{\bT}_i  & 0 \\
				0 & 0
			\end{bmatrix} + c_i^2\bI_4 - \trsp{\widehat{\bV}} \bQ_i \widehat{\bV}\bigg) \ba_i \geq 0, \ \forall \ba_i \in\bbR^4
		\end{split} \label{eq:4to3-aTa}
	\end{equation}
	In the above, we changed the coordinates by defining $\ba_i = \widehat{\bV} \bz_i\in\bbR^4$. The left-hand side of \eqref{eq:4to3-aTa} is a quadratic polynomial in entries of $\ba_i=[\balpha_i;a_{i4}]$, and we see in \eqref{eq:4to3-aTa} that $a_{i4}$ appears only in the term $\Big(c_i^2 - \trsp{\hat{\bw}_0} \bQ_i \hat{\bw}_0 \Big)a_{i4}^2 - 2a_{i4} \trsp{\hat{\bw}_0} \bQ_i \widehat{\bV}_0\balpha_i$. Since $c_i^2 > \trsp{\hat{\bw}_0} \bQ_i \hat{\bw}_0$, the left-hand side of \eqref{eq:4to3-aTa} is minimized as a function of $a_{i4}$ at $a_{i4}=\trsp{\hat{\bw}_0} \bQ_i \widehat{\bV}_0\balpha_i /\big(c_i^2 - \trsp{\hat{\bw}_0} \bQ_i \hat{\bw}_0 \big)$. Thus, with  $\bC_i := c_i^2\bI_3 -  \frac{\trsp{\widehat{\bV}_0}\bQ_i \hat{\bw}_0 \trsp{\hat{\bw}_0} \bQ_i \widehat{\bV}_0}{c_i^2- \trsp{\hat{\bw}_0} \bQ_i \hat{\bw}_0}$, we get
	\begin{equation*}
		\trsp{\ba_i} \bigg( \begin{bmatrix}
			\widehat{\bT}_i  & 0 \\
			0 & 0
		\end{bmatrix} + c_i^2\bI_4 - \trsp{\widehat{\bV}} \bQ_i \widehat{\bV}\bigg) \ba_i	\geq \trsp{\balpha_i} \Big(\widehat{\bT}_i + \bC_i - \trsp{\widehat{\bV}_0} \bQ_i \widehat{\bV}_0\Big) \balpha_i
	\end{equation*}	
	with the equality attained at $a_{i4}=\trsp{\hat{\bw}_0} \bQ_i \widehat{\bV}_0\balpha_i /\big(c_i^2 - \trsp{\hat{\bw}_0} \bQ_i \hat{\bw}_0 \big)$. This proves
	\begin{equation*}
		\widehat{\bS}_{i} + c_i^2\bI_4 -  \bQ_i \succ 0 \Leftrightarrow \widehat{\bT}_i + \bC_i - \trsp{\widehat{\bV}_0} \bQ_i \widehat{\bV}_0 \succ 0 
	\end{equation*}
	and we will prove the latter, which with the definition of $\widehat{\bT}_i$ \eqref{eq:Ti} can be written as 
	\begin{equation*}
		\bC_i - \frac{\Big( \sum_{i=1}^{\ell} \trsp{\hat{\bw}_0}\bQ_i\hat{\bw}_0 \Big) \bI_3}{\ell} + \frac{\sum_{j\neq i} \Big( \trsp{\widehat{\bV}_0} \bQ_{j} \widehat{\bV}_0 - \trsp{\widehat{\bV}_0} \bQ_{i} \widehat{\bV}_0 \Big)}{\zeta(\ell - 1)} \succ 0
	\end{equation*}
	Note that the above can be written as $\trsp{\widehat{\bV}_0} (\cdot) \widehat{\bV}_0 \succ 0$, so a sufficient condition for the above to hold is $(\cdot)\succ 0$; in other words, it suffices to prove
	\begin{equation*}
		c_i^2\bI_3 -  \frac{\bQ_i \hat{\bw}_0 \trsp{\hat{\bw}_0} \bQ_i }{c_i^2- \trsp{\hat{\bw}_0} \bQ_i \hat{\bw}_0} - \frac{\Big( \sum_{i=1}^{\ell} \trsp{\hat{\bw}_0}\bQ_i\hat{\bw}_0 \Big) \bI_3}{\ell} + \frac{\sum_{j\neq i} \big(  \bQ_{j}  -  \bQ_{i} \big)}{\zeta(\ell - 1)} \succ 0
	\end{equation*}
	Since $\norm{\bQ_i \hat{\bw}_0}^2\geq \lambda_{\textnormal{max}}(\bQ_i \hat{\bw}_0 \trsp{\hat{\bw}_0} \bQ_i)$, the above holds whenever we have
	\begin{equation*}
		c_i^2 -  \frac{ \norm{\bQ_i \hat{\bw}_0}^2 }{c_i^2- \trsp{\hat{\bw}_0} \bQ_i \hat{\bw}_0} - \frac{ \sum_{i=1}^{\ell} \trsp{\hat{\bw}_0}\bQ_i\hat{\bw}_0 }{\ell} - \frac{\lambda_{\textnormal{max}}\Big(\sum_{j\neq i} \big(  \bQ_{i}  -  \bQ_{j} \big)\Big)}{\zeta(\ell - 1)} > 0,
	\end{equation*}
	Using the definition of $d_i$ in \eqref{eq:noisy+outlier-free2-condition}, the above can be written as
	\begin{equation*}
		c_i^2- \trsp{\hat{\bw}_0} \bQ_i \hat{\bw}_0 - \frac{ \norm{\bQ_i \hat{\bw}_0}^2 }{c_i^2- \trsp{\hat{\bw}_0} \bQ_i \hat{\bw}_0} - d_i  > 0.
	\end{equation*}
	Multiplying the above inequality by $c_i^2- \trsp{\hat{\bw}_0} \bQ_i \hat{\bw}_0$ leads to a quadratic inequality in $c_i^2- \trsp{\hat{\bw}_0} \bQ_i \hat{\bw}_0$. Solve it for $c_i^2- \trsp{\hat{\bw}_0} \bQ_i \hat{\bw}_0$ and we get equivalently
	\begin{equation*}
		c_i^2 > \trsp{\hat{\bw}_0} \bQ_i \hat{\bw}_0 + \frac{ d_i  + \sqrt{d_i^2 + 4\norm{\bQ_i \hat{\bw}_0}^2}}{2},
	\end{equation*}
	which holds true, as long as our assumption \eqref{eq:noisy+outlier-free2-condition} holds. This finishes the proof.
\end{proof}

\section{Supplementary Material: Detailed Comparison to Yang \& Carlone \cite{Yang-ICCV2019}}\label{subsection:comparision-YC}
We mentioned that our \eqref{eq:SDP} can be treated as equivalent to the \textit{naive} relaxation of Yang \& Carlone (see Proposition 6 of  \cite{Yang-ICCV2019}). Here we review their relaxation, highlight the subtle differences, and conclude that our relaxation is tight if an only if so is theirs. 

Their derivation starts with \eqref{eq:TLS-R} and first arrives at (an equivalent version of) \eqref{eq:TLS-Q} via the use of quaternion products, which though could be  avoided in view of relation \eqref{eq:rotation-quaternion}. Then they defined $\theta'_i\in\{\pm 1\}$ and $\ba_i := \theta_i'\bw_0$ (so $\theta_i'=\trsp{\ba_i}\bw_0$), leading to the following equivalence (as we simplify and translate):
\begin{equation*}
	\begin{split}
		\eqref{eq:TLS-Q} \Leftrightarrow& \min_{\bw_0\in\bbS^3,\theta_i'\in\{\pm 1\}} \sum_{i=1}^{\ell}  \frac{1+\theta_i'}{2} \trsp{\bw}_0\bQ_i \bw_0 + \sum_{i=1}^{\ell}  \frac{1-\theta_i'}{2} c_i^2   \label{eq:Yang-1}\\
		\Leftrightarrow& \min_{\bw_0\in\bbS^3,\ba_i\in\{\pm \bw_0\}} \frac{\sum_{i=1}^{\ell}\trsp{\bw}_0\bQ_i \bw_0}{2} + \frac{\sum_{i=1}^{\ell}\trsp{\ba}_i\bQ_i \bw_0}{2} + \sum_{i=1}^{\ell}\frac{c_i^2}{2} (1-\theta_i')  \\
		\Leftrightarrow& \min_{\bw_0\in\bbS^3,\ba_i\in\{\pm \bw_0\}} \frac{\sum_{i=1}^{\ell}\trsp{\ba}_i\bQ_i \ba_i}{2} + \frac{\sum_{i=1}^{\ell}\trsp{\ba}_i\big(\bQ_i- c_i^2\bI_4 \big) \bw_0}{2} + \frac{\sum_{i=1}^{\ell}c_i^2}{2}
	\end{split}
\end{equation*}
Define $\pa:=[\bw_0;\ba_1;\dots;\ba_\ell]$, use our definition \eqref{eq:def-tildepQ} of $\pQ$, and let $\pQ'\in\bbR^{4(\ell+1)\times 4(\ell+1)}$ be a block diagonal matrix with $[\pQ']_{00}=0$ and $[\pQ']_{ii}=\bQ_i$. Then \eqref{eq:Yang-1} becomes
\begin{equation}\label{eq:YC-QCQP}
	\begin{split}
		\min_{\pa\in\bbR^{4(\ell+1)}} &\ \ \frac{1}{2} \trace\Big(\pQ' \pa \trsp{\pa}\Big) + \frac{1}{2}\trace\Big(\pQ \pa \trsp{\pa}\Big)  + \frac{ \sum_{i=1}^{\ell}c_i^2 }{2}  \\
		\text{s.t.}& \ \ [\pa\trsp{\pa}]_{00} = [\pa\trsp{\pa}]_{ii}, \ \ \forall\ i\in\{1,\dots,\ell\} \\
		& \ \ \trace\big([\pa\trsp{\pa}]_{00}\big) =  1 
	\end{split}\tag{QCQP-YC}
\end{equation}
Here, we used the notation $\pQ$ defined in \eqref{eq:def-tildepQ}. To summarize, \eqref{eq:YC-QCQP} is equivalent to the non-convex QCQP of Yang \& Carlone \cite{Yang-ICCV2019} (Proposition 4 \cite{Yang-ICCV2019}). Moreover, \eqref{eq:YC-QCQP} is also equivalent to our \eqref{eq:P-QCQP}, meaning that globally minimizing one of them solves the other and their minimum values $\hat{g}_{\text{\ref{eq:P-QCQP}}}$ and $\hat{g}_{\text{\ref{eq:YC-QCQP}}}$ are both equal to the minimum  $\hat{g}_{\text{\ref{eq:TLS-Q}}}$ of \eqref{eq:TLS-Q}. Relax \eqref{eq:YC-QCQP}, and we get
\begin{equation}\label{eq:SDP-YC} 
	\begin{split}
		\min_{\pA\succeq 0} &\ \ \frac{1}{2}\trace\big(\pQ' \pA\big) + \frac{1}{2}\trace\big(\pQ \pA\big) + \frac{ \sum_{i=1}^{\ell}c_i^2 }{2}  \\
		\text{s.t.}& \ \ [\pA]_{00} = [\pA]_{ii}, \ \ \forall\ i\in\{1,\dots,\ell\}  \\
		& \ \ \trace\big([\pA]_{00}\big) = 1
	\end{split} \tag{SDR-YC} 
\end{equation}
\eqref{eq:SDP-YC} is equivalent to the naive relaxation of Yang \& Carlone \cite{Yang-ICCV2019} in their Proposition 6. Not very obvious is the equivalence between \eqref{eq:SDP} and \eqref{eq:SDP-YC}:
\begin{proposition}\label{prop:SDP-equivalence}
	Assume $c_i^2\neq \lambda_{\textnormal{min}}(\bQ_i)$, $c_i^2\neq \lambda_{\textnormal{max}}(\bQ_i)$, $\forall i=1,\dots,\ell$. If \eqref{eq:TLS-Q} rejects all outliers and keeps all inliers, then \eqref{eq:SDP} is tight $\Leftrightarrow$  \eqref{eq:SDP-YC} is tight. 
\end{proposition}

\begin{proof}
	As before, assume that the first $k^*$ point pairs are inliers. We first prove the $\Leftarrow$ direction. Assume \eqref{eq:SDP-YC} is tight, i.e., it admits $\hat{\pa}\trsp{\hat{\pa}}$ as an global minimizer, where $\hat{\pa}$ minimizes \eqref{eq:YC-QCQP} (Definition \ref{def:tightness}). Since \eqref{eq:TLS-Q} rejects all outliers and keeps all inliers, $\hat{\pa}$ is of the form $[\hat{\bw}_0; \hat{\bw}_0; \dots; \hat{\bw}_0; -\hat{\bw}_0; \dots; -\hat{\bw}_0]$; here $\hat{\bw}_0$ is an eigenvector of $ \sum_{i=1}^{k^*} \bQ_i $ corresponding to its minimum eigenvalue $ \sum_{i=1}^{k^*} \trsp{\hat{\bw}_0}\bQ_i \hat{\bw}_0$ and $\hat{\bw}_0$ globally minimizes \eqref{eq:TLS-Q}. Moreover, we have $\hat{g}_{\text{\ref{eq:P-QCQP}}} = \big( \sum_{i=1}^{k^*} \trsp{\hat{\bw}_0}\bQ_i \hat{\bw}_0 \big) + \sum_{j=k^*+1}^{\ell} c_j^2$. Since \eqref{eq:SDP-YC} is tight, the dual program of \eqref{eq:SDP-YC}, or (derivation omitted)
	\begin{equation}\label{eq:Dual-YC}
		\max_{\mu,\pY}\ \ \mu +  \frac{ \sum_{i=1}^{\ell}c_i^2 }{2} \ \ \ \ \ \text{s.t.} \ \ \ \ \ \frac{1}{2}\pQ'+ \frac{1}{2}\pQ-\mu\pB-\pY \succeq 0 \tag{D-YC},
	\end{equation}
	has the minimum value $\sum_{i=1}^{k^*} \trsp{\hat{\bw}_0}\bQ_i \hat{\bw}_0 + \sum_{j=k^*+1}^{\ell} c_j^2$; here $\pB$ was defined in \cref{lemma:SDP-D} and $\pY$ is the $4(\ell+1)\times 4(\ell+1)$ block diagonal matrix satisfying 
	\begin{equation}\label{eq:Y-constraints}
		\begin{cases}
			\text{$\pY$ is symmetric}, \ \ [\pY]_{00} + \sum_{i=1}^{\ell} [\pY]_{ii} =0, \ \ \forall\ i\in\{1,\dots,\ell\} \\
			\text{all other entries of $\pY$ are zero.}
		\end{cases}
	\end{equation}
	In particular, the optimal $\hat{\pY}$ is of the form \eqref{eq:Y-constraints} that fulfills the optimality conditions:
	\begin{equation}\label{eq:KKT-YC}
		\begin{cases}
			\Big( \frac{1}{2}(\pQ'+\pQ)-\big(\sum_{i=1}^{k^*} \trsp{\hat{\bw}_0}\bQ_i \hat{\bw}_0 +\sum_{j=k^*+1}^{\ell} c_j^2  - \frac{ \sum_{i=1}^{\ell}c_i^2 }{2} \big)\pB-\hat{\pY} \Big) \hat{\pa} = 0\\		
			\frac{1}{2}(\pQ'+\pQ)-\big(\sum_{i=1}^{k^*} \trsp{\hat{\bw}_0}\bQ_i \hat{\bw}_0 +\sum_{j=k^*+1}^{\ell} c_j^2  - \frac{ \sum_{i=1}^{\ell}c_i^2 }{2} \big)\pB-\hat{\pY} \succeq 0
		\end{cases} 
	\end{equation}
	Similarly to \eqref{eq:stationary2} of \cref{prop:tightnesscondition}, the first condition of \eqref{eq:KKT-YC} is equivalent to 
	\begin{equation}\label{eq:stationary-noisy+outlier-YC}
		\begin{cases}
			\big( 3\bQ_i - c_i^2\bI_4 - 4[\hat{\pY}]_{ii} \big) \hat{\bw}_0 = 0, & i\in \{1,\dots, k^* \} \\
			\big( \bQ_j + c_j^2\bI_4 - 4[\hat{\pY}]_{jj} \big) \hat{\bw}_0 = 0, & j\in \{k^*+1,\dots, \ell \} 
		\end{cases}
	\end{equation}
	We will prove $\hat{\pw} \trsp{\hat{\pw}}$ globally minimizes \eqref{eq:SDP}; here $\hat{\pw}:=[\hat{\bw}_0;\dots;\hat{\bw}_0;0;\dots;0]$. Let  $[\hat{\pD}]_{0i}:=\bQ_i - 2[\hat{\pY}]_{ii}$ for every $i=1,\dots,\ell$ and $\hat{\mu}:=\sum_{i=1}^{k^*} \trsp{\hat{\bw}_0}\bQ_i \hat{\bw}_0-\sum_{i=1}^{k^*}c_i^2$. It now suffices to show that the stationary condition \eqref{eq:stationary2} and dual feasibility condition \eqref{eq:dualfeasibility2} of \cref{prop:tightnesscondition} are fulfilled. Note that \eqref{eq:stationary2} follows directly from \eqref{eq:stationary-noisy+outlier-YC} and the definition of $[\hat{\pD}]_{0i}$. It remains to prove that the second condition of \eqref{eq:KKT-YC} implies \eqref{eq:dualfeasibility2}. First observe that the second condition of \eqref{eq:KKT-YC}, multiplied by $4$, implies ($\forall \bz_i\in\bbR^4$)
	\begin{equation}\label{eq:KKT-YC2}
		\begin{cases}
			\trsp{\bz_0}\big( -4\hat{\mu}\bI_4- 2\sum_{i=1}^{\ell}c_i^2\bI_4 - 4[\hat{\pY}]_{00} \big) \bz_0 + \\
			\sum_{i=1}^{\ell} \trsp{\bz_i}\big( 2\bQ_i - 4[\hat{\pY}]_{ii} \big) \bz_i - \\
			\sum_{i=1}^{\ell} 2\trsp{\bz_0}\big( c_i^2\bI_4 -  \bQ_i  \big) \bz_i
		\end{cases} \geq 0
	\end{equation}
	Similarly to \cref{lemma:D0i-PSD}, since $c_i^2\neq \lambda_{\textnormal{min}}(\bQ_i)$ and $c_i^2\neq \lambda_{\textnormal{max}}(\bQ_i)$, we have $2[\hat{\pD}]_{0i}:=2\bQ_i - 4[\hat{\pY}]_{ii}\succ 0$. now, substitute the minimizer  $\bz_i=(2[\hat{\pD}]_{0i})^{-1}\big( c_i^2\bI_4 -  \bQ_i  \big)\bz_0$ into \eqref{eq:KKT-YC2}, and we see that \eqref{eq:KKT-YC2} divided by $4$ is equivalent to ($\forall \bz_0\in\bbR^4$)
	\begin{equation*}
		\trsp{\bz_0}\Big(-\hat{\mu}\bI_4- \sum_{i=1}^{\ell}\frac{c_i^2\bI_4}{2} - [\hat{\pY}]_{00} - \sum_{i=1}^{\ell} \frac{\big( c_i^2\bI_4 -  \bQ_i  \big) (2[\hat{\pD}]_{0i})^{-1} \big( c_i^2\bI_4 -  \bQ_i  \big)}{4}  \Big)\bz_0 \geq 0,
	\end{equation*}
	or equivalently (with $\bz_0$ now removed)
	\begin{equation*}
		\begin{split}
			&-\hat{\mu}\bI_4- \sum_{i=1}^{\ell}\frac{c_i^2\bI_4}{2} - [\hat{\pY}]_{00} - \sum_{i=1}^{\ell} \frac{\big( c_i^2\bI_4 -  \bQ_i  \big) (2[\hat{\pD}]_{0i})^{-1} \big( c_i^2\bI_4 -  \bQ_i  \big)}{4}  \succeq 0 \\
			\Leftrightarrow&-\hat{\mu}\bI_4- \sum_{i=1}^{\ell}\bigg( \frac{c_i^2\bI_4}{2} - [\hat{\pY}]_{ii} +  \frac{\big( c_i^2\bI_4 -  \bQ_i  \big) (2[\hat{\pD}]_{0i})^{-1} \big( c_i^2\bI_4 -  \bQ_i  \big)}{4}  \bigg) \succeq 0 \\
			\Leftrightarrow& -\hat{\mu}\bI_4- \sum_{i=1}^{\ell}  \frac{ 2c_i^2\bI_4 - 4[\hat{\pY}]_{ii} + \big( c_i^2\bI_4 -  \bQ_i  \big) (2[\hat{\pD}]_{0i})^{-1} \big( c_i^2\bI_4 -  \bQ_i  \big)}{4}  \succeq 0 \\
			\Leftrightarrow& -\hat{\mu}\bI_4- \sum_{i=1}^{\ell}  \frac{   2[\hat{\pD}]_{0i} + 2c_i^2\bI_4 - 2\bQ_i + \big( c_i^2\bI_4 -  \bQ_i  \big) (2[\hat{\pD}]_{0i})^{-1} \big( c_i^2\bI_4 -  \bQ_i  \big)}{4}  \succeq 0 \\
			\Leftrightarrow& -\hat{\mu}\bI_4 - \sum_{i=1}^{\ell} \frac{\big(2[\hat{\pD}]_{0i} - \bQ_i+c_i^2\bI_4\big)\big(2[\hat{\pD}]_{0i}\big)^{-1}\big(2[\hat{\pD}]_{0i} - \bQ_i+c_i^2\bI_4\big)}{4} \succeq 0,
		\end{split}
	\end{equation*}
	or equivalently (with $\bz_0$ back)
	\begin{equation*}
		\trsp{\bz_0}\bigg( -\hat{\mu}\bI_4 - \sum_{i=1}^{\ell} \frac{\big(2[\hat{\pD}]_{0i} - \bQ_i+c_i^2\bI_4\big)\big(2[\hat{\pD}]_{0i}\big)^{-1}\big(2[\hat{\pD}]_{0i} - \bQ_i+c_i^2\bI_4\big)}{4}\bigg)\bz_0 \geq 0
	\end{equation*}
	or equivalently (with all other $\bz_i$'s back)
	\begin{equation}
		-\hat{\mu} \norm{\bz_0}^2 + 2\sum_{i=1}^{\ell}\trsp{\bz_i}[\hat{\pD}]_{0i} \bz_i -\sum_{i=1}^{\ell}\trsp{\bz_0}\big(2[\hat{\pD}]_{0i} - \bQ_i+c_i^2\bI_4\big)\bz_i\geq 0,
	\end{equation}
	which is exactly \eqref{eq:dualfeasibility2}. The above derivation also implies that the $\Rightarrow$ direction is true.
\end{proof}
\begin{lemma}\label{lemma:D0i-PSD}
	Condition \eqref{eq:dualfeasibility2} implies $[\hat{\pD}]_{0i} \succeq 0$ for every $i=1,\dots,\ell$. Moreover, if $c_i^2\neq \lambda_{\textnormal{min}}\big(\bQ_i\big)$ and $c_i^2\neq \lambda_{\textnormal{max}}\big(\bQ_i\big)$ then \eqref{eq:dualfeasibility2} implies $[\hat{\pD}]_{0i}\succ 0$.
\end{lemma}
\begin{proof}[\cref{lemma:D0i-PSD}]
	Substitute $\bz_0=\bz_1=\cdots=\bz_{\ell-1}=0$ into \eqref{eq:dualfeasibility2} and we have $\trsp{\bz_\ell}[\hat{\pD}]_{0\ell} \bz_\ell\geq 0$. This proves that  \eqref{eq:dualfeasibility2}  implies $[\hat{\pD}]_{0i} \succeq 0$ for every $i=1,\dots,\ell$. 
	
	Assume If $c_1^2\neq \lambda_{\textnormal{min}}\big(\bQ_1\big)$ and $c_1^2\neq \lambda_{\textnormal{max}}\big(\bQ_1\big)$ and $[\hat{\pD}]_{01}$ has an eigenvalue $0$ with $[\hat{\pD}]_{01} \bz_1'=0$ for some $\bz_1'\in\bbS^3$. With some $t\in\bbR$, substitute $\bz_0\in\bbS^3$ and $\bz_1=t \bz_1'$ and $\bz_2=\cdots=\bz_\ell=0$ into \eqref{eq:dualfeasibility2}, and we get $-\hat{\mu} + t \trsp{\bz_0} (c_1^2\bI_4 - \bQ_1)\bz_1'\geq 0$ for every $\bz_0\in\bbS^3$. Since $c_1^2$ is not equal to any eigenvalue of $\bQ_1$ (\cref{lemma:rotation-quaternion}), the matrix $c_1^2\bI_4 - \bQ_1$ is of full rank, and we have $(c_1^2\bI_4 - \bQ_1)\bz_1'\neq 0$. Choose $\bz_0$ such that $\trsp{\bz_0} (c_1^2\bI_4 - \bQ_1)\bz_1'<0$, and choose $t\to \infty$, then we see that $-\hat{\mu} + t \trsp{\bz_0} (c_1^2\bI_4 - \bQ_1)\bz_1'< 0$, a contradiction. This proves that \eqref{eq:dualfeasibility2} and $c_i^2\neq \lambda_{\textnormal{min}}\big(\bQ_i\big)$ and $c_i^2\neq \lambda_{\textnormal{max}}\big(\bQ_i\big)$ imply $[\hat{\pD}]_{0i}\succ 0$.
\end{proof}

\bibliographystyle{siamplain}
\bibliography{/home/liangzu/Dropbox/Liangzu}

\end{document}